\documentclass{amsart}
\pdfoutput=1

\usepackage{graphicx}
\usepackage{mathptmx}    
\usepackage{latexsym}
\usepackage{amsfonts}
\usepackage{amsmath}

\newtheorem{theorem}{Theorem}[section]
\newtheorem{lemma}[theorem]{Lemma}
\newtheorem{corollary}[theorem]{Corollary}
\newtheorem{proposition}[theorem]{Proposition}

\theoremstyle{definition}

\newtheorem{example}[theorem]{Example}

\theoremstyle{remark}
\newtheorem{remark}[theorem]{Remark}

\numberwithin{equation}{section}

\begin{document}
\title{A geodesic stratification of two-dimensional semi-algebraic sets}
\author{Chengcheng Yang }
\address{Department of Mathematics, Rice University, 6100 Main Street, Houston TX 77005}
\email{cy2@rice.edu}           
\thanks{The author was support in part by NSF Award\#1745670}
\date{\today}
\dedicatory{This paper is dedicated to my parents and my thesis advisor Dr. Robert Hardt. My parents have supported me unconditionally. My advisor has given me many valuable advices.}

\begin{abstract}
Given any arbitrary semi-algebraic set $X$, any two points in $X$ may be joined by a piecewise $C^2$ path $\gamma$ of shortest length. Suppose $\mathcal{A}$ is a semi-algebraic stratification of $X$ such that each component of $\gamma \cap \mathcal{A}$ is either a singleton or a real analytic geodesic segment in $\mathcal{A}$, the question is whether $\gamma \cap \mathcal{A}$ has at most finitely many such components. This paper gives a semi-algebraic stratification, in particular a cell decomposition, of a real semi-algebraic set in the plane whose open cells have this finiteness property. This provides insights for high dimensional stratifications of semi-algebraic sets in connection with geodesics. 
\end{abstract}

\maketitle
\section{introduction}
A semi-algebraic set $X$ in the plane can be described as: 
\begin{equation*}
X = \bigcup_{i=1}^{I} \bigcap_{j=1}^{J} \{(x, y) \in \mathbb{R}^2: f_{i, j}(x,y) = 0, g_{i,j}(x,y ) > 0\}, 
\end{equation*}
where $f_{i,j}$, $g_{i,j}$ are polynomials in two variables. We see that $X$ is a finite union of sets in the form obtained by taking the intersection of an algebraic set (i.e. $\{f_1(x, y)= 0, \ldots, f_k(x,y)=0\}$) with an open semi-algebraic set (i.e. $\{g_1(x, y) > 0, \ldots, g_m(x, y)>0\}$).

The triangulability question for algebraic sets was first considered by van de Waerden in 1929  \cite{V}. It is a well-known theorem that every algebraic set is triangulizable \cite{H}. On the other hand, in 1957 Whitney introduced another splitting process that divides a real algebraic $V$ into a finite union of ``partial algebraic manifolds" \cite{W}. An algebraic partial manifold $M$ is a point set, associated with a number $\rho$, with the following property. Take any $p \in M$. Then there exists a set of polynomials $f_1, \ldots, f_{\rho}$, of rank $\rho$ at $p$, and a neighborhood $U$ of $p$, such that $M \cap U$ is the set of zeros in $U$ of these $f_i$. The splitting process uses the rank of a set $S$ of functions $f_1, \ldots, f_s$ at a point $p$, where the rank of $S$ at $p$ is the number of linearly independent differentials $df_1(p), \ldots, df_s(p)$. 

In 1975 Hironaka reproved that every semi-algebraic set is also triangulable and also generalized it to sub-analytic sets \cite{H}. His proof came from a paper of Lojasiewicz in 1964, in which Lojasiewicz proved that a semi-analytic set admits a semi-analytic triangulation \cite{L}. In 1975, Hardt also proved the triangulation result for sub-analytic sets by inventing another new method \cite{H2}. Since any semi-algebraic set is also semi-analytic, thus is sub-analytic, both Hironaka and Hardt's results showed that a semi-algebraic set admits a triangulation, that is to say, it is homeomorphic to the polyhedron of some simplicial complex. 

Following the examples of Whitney's stratification and Lojasiewicz/Hironaka/Hardt's triangulation, this paper tries to build a cell-complex stratification such that it admits an analytical condition concerning shortest-length curves. The idea is explained more precisely as follows. 

Suppose $A$, $B$ are two arbitrary points in $X$, and $\gamma$ is a piecewise $C^2$ curve from $A$ to $B$ lying entirely in $X$ such that its length is the shortest among all possible such curves. We search for a semi-algebraic cell decomposition (that is each cell is a semi-algebraic set in $\mathbb{R}^2$) $\mathcal{A}$ of $X$ , such that the intersection of $\gamma$ with every cell in $\mathcal{A}$ is either empty or consists of finitely many components, each of which is either a singleton or a geodesic line segment. We will explain the meaning of a geodesic line segment soon. If a cell decomposition $\mathcal{A}$ satisfies this property, we will simply say that $\mathcal{A}$ satisfies {\bf the finiteness property}. 

Our first step is to assume that $X$ is an arbitrary affine algebraic variety, that is, 
\begin{equation}\label{eqn:F0.1}
X = \{f_1(x, y)=0, \ldots, f_k(x,y)=0\}, k \geq 1,
\end{equation}
where the $f_i(x,y)$ are nonzero distinct polynomials. We argue that a cell decomposition $\mathcal{A}$ exists with the finiteness property, and $(X, \mathcal{A})$ can be shown to be a CW complex. More generally, we may assume that $X$ is a finite union of sets in the above form, then the same conclusion holds. 

Our next step is to look at an open planar semi-algebraic set $X$ in the form of:
\begin{equation}\label{eqn:F0.2}
X = \{g_1(x, y) > 0, \ldots, g_m(x, y) > 0\}, m \geq 1,
\end{equation}
where the $g_j(x, y)$ are nonzero distinct polynomials. More generally, we may assume that $X$ is a finite union of sets in the above form. 
A cell decomposition for such an $X$ can also be established with the desired finiteness property which is also a CW decomposition. 

Our third step is to search for a CW decomposition with the desired finiteness property for $X$, which is a finite union of sets in the form of:
\begin{equation}\label{eqn:F0.3}
\{f_1(x, y)=0, \ldots, f_k(x,y)=0, g_1(x, y) > 0, \ldots, g_m(x, y) > 0\}, k \geq 1, m \geq 1. 
\end{equation}

Our final step is to take a finite union of sets in the previous steps:
(\ref{eqn:F0.1}) and (\ref{eqn:F0.2});
(\ref{eqn:F0.1}) and (\ref{eqn:F0.3});
(\ref{eqn:F0.2}) and (\ref{eqn:F0.3});
(\ref{eqn:F0.1}) and (\ref{eqn:F0.2}) and (\ref{eqn:F0.3}).

\section{The stratification of an affine algebraic set in $\mathbb{R}^2$}
Suppose $X =\{(x,y) \in \mathbb{R}^2: f(x, y) = 0\}$ is an affine algebraic set in the plane.  Suppose $p=(x_0, y_0)$ is a nonsingular point of $X$, that is $df(p) \neq 0$. Without loss of generality assuming that $\frac{\partial f}{\partial y} (p) \neq 0$,  the implicit function theorem implies that there exist open intervals $I$, $J$ of $x_0$, $y_0$, respectively, and a differentiable function $g: I \rightarrow J$ such that $g(x_0)=y_0$ and 
\begin{equation*}
\{(x, y) \in I \times J \ | \ f(x, y) = 0\} = \{(x, g(x)) \in \mathbb{R}^2 \ | \ x \in I\}
\end{equation*}
\cite{RefM&R}. So we obtain a smooth parametrization $g$ for $X$ in an open neighborhood of $p$. In the paper \cite{Y}, we've shown how to construct a cell decomposition with the desired finiteness property in the closed region below the graph of $g$ under the assumption that $g$ is a polynomial function. More generally, the closed region could be replaced by an open region below the graph, and the polynomial function could be replaced by a smooth function with finitely many strict inflection and local minimum points. The following lemma verifies that g is in fact a real analytic (thus smooth) function over the open interval $I$. 

\begin{lemma}\label{lem:F2.1}
Suppose $g$ is the differentiable function given as before by the implicit function theorem for the polynomial function $f(x, y)=0$ at $p=(x_0, y_0)$, where $\frac{\partial f}{\partial y}(p) \neq 0$, then g is a real analytic function over the open interval $I$. 
\end{lemma}

\begin{proof}
Since $f(x, y)$ is a polynomial in two variables, we can consider the complex polynomial function $f(z, w)$, where $z$, $w$ are variables in $\mathbb{C}$. It follows that $f(z, w)$ is a holomorphic function in two variables. Now we can apply the holomorphic implicit function theorem, since $\frac{\partial f}{\partial w} \neq 0$ at $(x_0, y_0) \in \mathbb{C} \times \mathbb{C}$, the equation $f(z, w)=0$ has a unique holomorphic solution $w(z)$ in a neighborhood $x_0$ that satisfies $w(x_0)=y_0$ \cite{K}. Hence $g(x) = w(x)$ for $x \in I$ when $I$ is in this neighborhood. 
\end{proof}

\begin{corollary}\label{cor:F2.3}
If $[a, b]$ is a closed and bounded interval contained in $I$ and suppose $g$ is not a linear function, then $g$ has finitely many inflection and local minimum points over $[a, b]$. 
\end{corollary}

\begin{proof}
For local minimum points (more generally, critical points), differentiating both sides of the equation $f(x, g(x)) = 0$ yields:
\begin{equation}\label{eqn:F2.1}
f_x(x, g(x)) + f_y(x, g(x)) \cdot g'(x) = 0,
\end{equation}
where we use $f_x, f_y$ as short-hand notations for $\frac{\partial f}{\partial x}$, $\frac{\partial f}{\partial y}$, respectively. It implies that $g'(x) = 0$ if and only if $f_x(x, g(x)) = 0$, because we may assume that $f_y(x, g(x)) \neq 0$ for all $x \in I$ by continuity. 

Since $f_x$ is also a polynomial, $f_x(x, g(x))$ is a real analytic function over the interval $I$. Therefore $f_x(x, g(x))$ has isolated zeros unless it is identically equal to zero in which case $g$ is a constant function. Since $[a, b]$ is a compact interval, there are at most finitely many zeros of $f_x(x, g(x))$ over $[a, b]$, thus there are at most finitely many local minimum points (or critical points) of $g$ over $[a, b]$. 

Similarly, for inflection points, we differentiate equation~(\ref{eqn:F2.1}) again to obtain the following equation: 
\begin{equation}\label{eqn:F2.6}
f_{xx}(x, g(x)) + 2f_{xy}(x, g(x)) \cdot g'(x) + f_{yy}(x, g(x)) \cdot g'(x)^2 + f_y(x, g(x)) \cdot g''(x) = 0.
\end{equation}
It follows that $g''(x) = 0$ if and only if 
\begin{equation*}
f_{xx}(x, g(x)) + 2f_{xy}(x, g(x)) \cdot g'(x) + f_{yy}(x, g(x)) \cdot g'(x)^2 =0,
\end{equation*}
which is real analytic and so has at most finitely many zeros over any compact interval $[a, b]$ because $g''$ is not identically zero by hypothesis. 
\end{proof}

Now we are ready to give a cell decomposition for a compact and connected algebraic variety of an irreducible polynomial in two real variables. 
\begin{theorem}\label{thm:F2.1}
Suppose $f(x,y)$ is an irreducible polynomial function and $X = V(f)$ is the affine algebraic variety determined by the zeros of $f$. Assume $X$ is compact and connected, then $X$ has a cell decomposition $\mathcal{A}$ with the desired finiteness property. 
\end{theorem}

The proof immediately gives the following corollary.

\begin{corollary}\label{cor:F2.1}
Given the cell decomposition $\mathcal{A}$ as in the theorem, if $\gamma$ is any shortest-length piecewise $C^2$-curve between two points in $X$, then the intersection of $\gamma$ with any 0-cell in $\mathcal{A}$ is either empty or a singleton; the intersection of $\gamma$ with any 1-cell in $\mathcal{A}$ is  either empty, or a continuous line segment contained in the 1-cell.
\end{corollary}

\begin{remark}
A continuous line segment contained in a 1-cell is one example of a geodesic line segment. In general, a line segment is said to be {\bf geodesic} if it is either a straight line segment or a continuous line segment contained (partially or entirely) inside a 1-cell. 
\end{remark}

\begin{proof}
If $f$ is a polynomial in $x$- (or $y$-)variable only, then $X$ has at most one zero and the theorem follows trivially. Without loss of generality, we may assume that $f$ is a polynomial function that has both $x$ and $y$ variables. In particular, $f_x$ and $f_y$ are not zero. Then we have that $f$ and $f_x$ share no common factors, because $f$ is irreducible by hypothesis. An algebraic geometry theorem says if $k$ is an arbitrary commutative field, and $F, G \in k[x, y]$ are nonzero polynomials without common factors, then $V(F) \cap V(G) $ is finite \cite{P}. Applying this theorem, we conclude that 
\begin{equation}\label{eqn:F2.2}
V(f) \cap V(f_x) \text{ is finite.}
\end{equation}
Similarly, the same reasoning implies that
\begin{equation}\label{eqn:F2.3}
V(f) \cap V(f_y) \text{ is finite.}
\end{equation}
In particular, the set $S$ of singular points in $X$, that is 
\begin{equation*}
S = \{ x \in X: \text{ where }\frac{\partial f}{\partial x} (x) = 0  \text{ and } \frac{\partial f}{\partial y} (x)= 0\}
\end{equation*}
consists of at most finitely many points. 

{\it Case 1:} suppose $S$ is an empty set, then at every $p= (x_0, y_0) \in X$, either $f_x \neq 0$ or $f_y \neq 0$. We've known from Lemma~\ref{lem:F2.1} that $p$ has an open neighborhood $I \times J$ whose intersection with $X$ is the graph of a real analytic function over either $I$ or $J$. Shrinking $I$ and $J$ if necessary, we may also assume that the intersection of the closed neighborhood $\bar{I} \times \bar{J}$ with $X$ is the graph of a real analytic function.

Since $X$ is compact, $X$ can be covered by finitely many such open neighborhoods, say $I_1 \times J_1$, \ldots, $I_r \times J_r$, where $I_i \times J_i \not\subset I_j \times J_j$ for $i \neq j$. In each intersection of $X$ with $\bar{I_i} \times \bar{J_i}$, the graph has finitely many (strict) inflection and local minimum points, thus giving a cell decomposition as we've shown in [put a book citation here]. More precisely in this special case, the 0-cells are the (strict) inflection and local minimum points, together with the two endpoints; and the 1-cells are the graphs in between them. Therefore, we find a finite cell decomposition for $X$ with the desired finiteness property. 

{\it Case 2:} suppose $S$ is not an empty set,  we know that $S$ consists of finitely many points, thus each of which is an isolated point in $X$. Let's pick an open ball $B(q)$ for each $q \in S$ such that $B(q)$ contains no other point in $S$. Furthermore, in virtue of (\ref{eqn:F2.2}) and (\ref{eqn:F2.3}) and shrinking $B(q)$ if necessary, we may assume that for every point $p$ in the closed ball $\bar{B}(q)$, $p \neq q$, we have
\begin{equation}
\label{eqn:F2.4}
\frac{\partial f}{\partial x}(p) \neq 0, \text{ and }  \ \frac{\partial f}{\partial x}(p) \neq 0.
\end{equation}
Again by the compactness of $X$, $X$ can be covered by finitely many open sets in the form of either $I \times J$ around a non-singular point $p$ or $B(q)$ for a singular point $q$ in $S$. For the intersection of $X$ with $\bar{I} \times \bar{J}$, we use the same cell decomposition as shown in case 1 above. For the intersection of $X$ with $B(q)$, we need to first prove the following lemma.

\begin{lemma}\label{lem:F2.2}
Under the same assumption as before, the intersection of the punctured ball $B(q) \setminus \{q\}$, for each $q \in S$, with $X$ consists of finitely many connected components, each of which is homeomorphic to the real line $\mathbb{R}$. 
\end{lemma}

\begin{proof}
Since $B(q) \setminus \{q\}$ is an open subset of $\mathbb{R}^2$, its intersection with $X$ is an open subset of $X$, thus consisting of open connected components. Each connected component is locally Euclidean due to (\ref{eqn:F2.4}). Using the subspace topology inherited from $\mathbb{R}^2$, each connected component is also Hausdorff and second-countable. Therefore, each connected component is a connected 1-dimensional manifold. By the classification theorem for connected 1-manifolds, each connected component is homeomorphic to $\mathbb{S}^1$ if it is compact, and $\mathbb{R}$ if it is not \cite{L2}. Suppose there exists a connected component $P$ that is homeomorphic to $\mathbb{S}^1$, then $P$ is compact in $X$, thus closed in $X$ (by the Hausdorff property of $X$). Because $P$ is also open in $X$, it follows that $P$ is both open and closed in $X$. By hypothesis $X$ is connected, so $X = P$. However, $q$ is also in $X$ and $q$ is not contained in $P$, so $X \neq P$, which is a contradiction. So every connected component in the intersection of $B(q) \setminus \{q\}$ is homeomorphic to the real line $\mathbb{R}$. 

Next we want to show that there are finitely many such components. Given a connected component $P$, $P$ is contained in the punctured open ball $B(q) \setminus \{q\}$, so its boundary is contained in the boundary of $B(q) \setminus \{q\}$, which is $\partial B(q) \cup \{q\}$. There are four possibilities for the two endpoints of $P$: they are on the circle $\partial B(q)$ and are the same; they are on the circle but not the same; one of them is on the circle and the other is $q$; they are both equal to $q$ (see Figure~\ref{fig:F2.1}). 

\begin{figure}[ht]
\includegraphics[width=13cm]{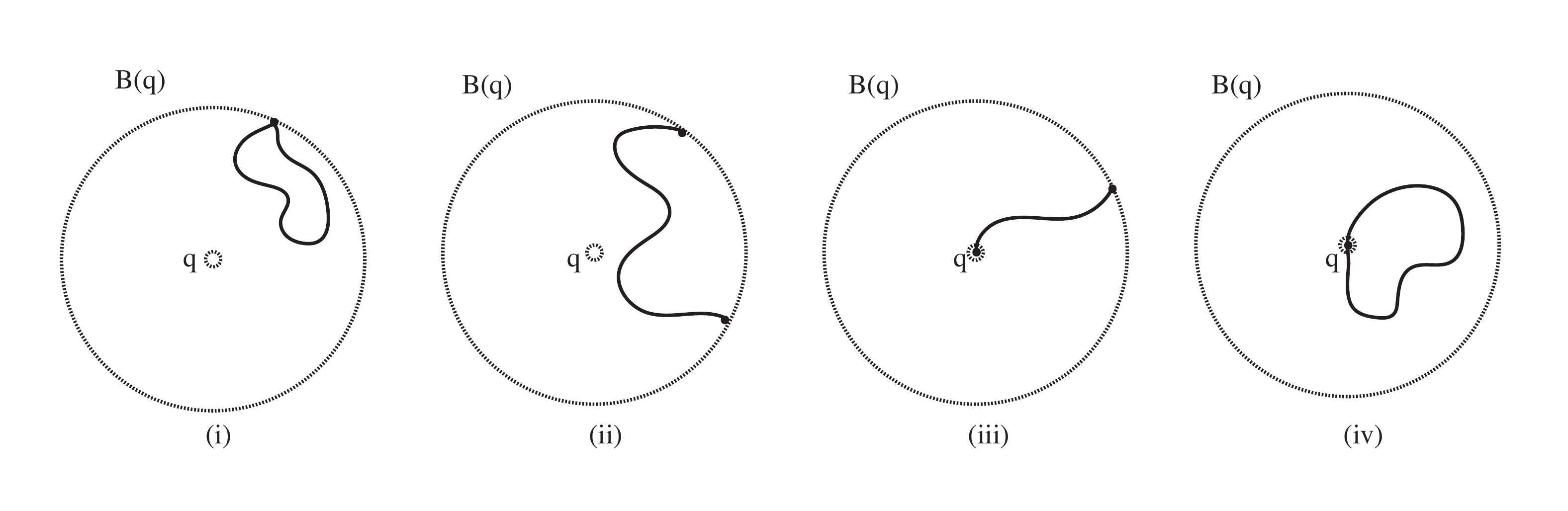}
\caption{Four possibilities of a connected component in $B(q) \setminus \{q\} \cap X$.}
\label{fig:F2.1}
\end{figure}

(i) Let's start with showing that the first situation is impossible. If the two endpoints of $P$ are both equal to a point $p$ on the boundary of $B(q)$, according to (~\ref{eqn:F2.4}), the closure $\bar{P}$ of $P$ is locally Euclidean, thus is homeomorphic to $\mathbb{S}^1$, because it is also closed and bounded. $\bar{P}$ is closed in $X$. Furthermore, there exists an open neighborhood $I \times J$ around $p$ such that the intersection of $X$ with $I \times J$ is equal to $\bar{P} \cap X$. It follows that $\bar{P}$ is also an open subset of $X$. Therefore, $X = \bar{P}$, since $X$ is connected by hypothesis, which yields a contradiction since $q \notin \bar{P}$. 

(ii) We show that there are at most finitely many components whose two endpoints are on the boundary of $B(q)$ which are not the same. Suppose for the sake of contradiction that there are infinitely many such components, then their endpoints are infinitely many, because each point on $\partial B(q)$ is an endpoint of at most two connected components. By sequential compactness of $\partial B(q)$, there exists a subsequence $\{p_j\}_{j=1}^{\infty}$ of these endpoints that converges to a point $p$ in $\partial B(q)$, which is also in $X$ for $X$ is closed. Since $p$ is non-singular, there exists an open neighborhood $I \times J$ of $p$ such that the intersection of $I \times J$ with $X$ is the graph of a real analytic function $g$. Then $g$ intersects the circle $\partial B(q)$ infinitely many times near $p$. Because $(1-t)^{1/2}$ has a power series expansion for $|t| < 1$, and it converges absolutely and uniformly on compact subsets of $(-1, 1)$, thus $(1-t)^{1/2}$ is real analytic over the open interval $(-1,1)$ ([put a book citation here, Folland, exercise 66, p. 139]). It follows that the arc near $p$ can be parametrized by a real analytic function as well. Since their difference is also real analytic and they have a zero that is not isolated, the graph of $g$ coincides with the circle near $p$. This implies that near $p$ there cannot be any endpoint of a connected component, thus leading to a contradiction. 

(iii) Similarly, there are at most finitely many components whose endpoints are made of one point on the boundary of $B(q)$ and one point being $q$. 

(iv) We finish the proof of the lemma by showing that the case when the two endpoints of $P$ are both equal to $q$ does not happen as well. Since $\frac{\partial f}{\partial x}$ and $\frac{\partial f}{\partial y}$ are both nonzero for every point in $B(q) \setminus \{q\}$, each point in $P$ has an open neighborhood $I \times J$ such that the intersection of $X$ with $I \times J$ is not only the graph of a real analytic function $g(x)$ for $x \in I$, but also the graph of a real analytic function $h(y)$ for $y \in J$. It follows that for each $x \in I$, $x$ satisfies the equation $x = h(g(x))$, thus implying  
\begin{equation}\label{eqn:F2.5}
1 = h'(y) \cdot g'(x), \text{ for } y = h(x), \text{ and } x \in I. 
\end{equation}
So $g'(x)$ is either $> 0$ or $< 0$ over the entire interval $I$. Without loss of generality, let us assume that $g' > 0$ for at least one point in $P$. Then consider the set $S$ of all points in $P$ satisfying the same property. That is, 
\begin{eqnarray}
S = \big\{r \in P: r \text{ has an open neighborhood } I \times J \text{ such that }  X \cap I \times J \text{ is the graph of}  \nonumber \\
 \text{ some real analytic function } g(x), \text { and } g'(x) > 0  \text{ for } x \in I\big\}. 
 \nonumber
\end{eqnarray}
$P$ being path-connected implies that $S = P$, because it is easy to see that $S$ is both open and closed, and $S$ is also non-empty. Since $P$ is homeomorphic to $\mathbb{R}$, we can choose an orientation for $P$ where locally the graph is increasing as we move along this direction. Thus, starting from $q$ and following this orientation, the y-coordinate always increases, therefore it is impossible that the other endpoint of $P$ returns to $q$. 

Finally, since only cases (ii) and (iii) are allowed, and there can be at most finitely many connected components in each case, we prove the lemma. 
\end{proof}

Let's continue proving the theorem. According to Lemma~\ref{lem:F2.2}, since there are at most finitely many connected components in the intersection of $X$ with $B(q) \setminus \{q\}$, we may shrink the open ball $B(q)$ if necessary to make sure that no component in case (ii) appears. Therefore, it remains to describe a cell decomposition for each component $P$ in case (iii) of the lemma.

Before proceeding with the description, we need to demonstrate the following lemma.
\begin{lemma}\label{lem:F2.3}
Suppose $P$ is a connected component in the intersection of $X$ with $B(q) \setminus \{q\}$, and one endpoint of $P$ is equal to $q$, and the other is on the boundary of $B(q)$, then $P$ is the graph of a real analytic function $g(x)$ for $x$ over an open interval $I= (a, b)$, where $-\infty < a < b < \infty$. Moreover, $g(x)$ has no critical points and at most finitely many strict inflection points over the interval $I$. 
\end{lemma}  

\begin{proof}
For every point $r$ on $P$, $r$ has an open neighborhood $I_r \times J_r$ such that the intersection of $X$ with $I_r \times J_r$ is the graph of a real analytic function $g_r(x)$ for $x \in I_r$. Moreover, we may assume that $I_r \times J_r$ is contained inside the open ball $B(q)$. Let 
\begin{equation*}
I = \bigcup_{r \in P} \  I_r.
\end{equation*}
Then we show that $I$ is connected, thus $I$ is an open interval. Suppose not, there exist two disjoint subsets $P_1$, $P_2$ of $P$, such that $P_1 \cup P_2 = P$, and 
\begin{equation*}
I = \big(\bigcup_{r \in P_1} \ I_r\big) \ \  \bigcup \ \  \big(\bigcup_{r \in P_2} \ I_r\big), 
\end{equation*}
where the union is disjoint. We can deduce that 
\begin{equation*}
\bigcup_{r \in P_1} \ I_r \times J_r \text{ and } \bigcup_{r \in P_2} \ I_r \times J_r \text{ are disjoint.}
\end{equation*}
Since 
\begin{equation*}
P = P \cap \bigcup_{r \in P} \ I_r \times J_r = \big(P \cap \bigcup_{r \in P_1} \ I_r \times J_r\big) \ \ \bigcup \ \ \big( P \cap \bigcup_{r \in P_2} \ I_r \times J_r \big),
\end{equation*}
$P$ turns out to be a disconnected set, which is a contradiction to the hypothesis that $P$ is a connected component. Thus $I$ is connected. Since $I$ is also open and bounded (under the extra assumption that each $I_r \times J_r$ is contained in $B(q)$), there exist two real numbers $a < b$ such that $I = (a, b)$. 

Given two distinct points $r_1, r_2$ in $P$, suppose that $I_{r_1}$ and $I_{r_2}$ overlap nontrivially, then we claim that $g_{r_1}$ and $g_{r_2}$ agree over the intersection $I_{r_1} \cap I_{r_2}$. 
The proof of the claim is essentially the same as shown in (\ref{eqn:F2.5}), except that in this case we look at the $x$-coordinate instead of the $y$-coordinate. As a result, $P$ can have one and only one graph over each point in the common interval of $I_{r_1}$ and $I_{r_2}$. It follows that $g_{r_1} = g_{r_2}$ there. 

Given this, for each $x \in I$, if we define the value of $g(x)$ to be $g_r(x)$ whenever $x \in I_r$ for some $r \in P$. Then $g$ is well-defined over the entire open interval $I$. Furthermore, $g$ is a real analytic function. 

(The following was first suggested to me by Dr. Hardt, which significantly simplifies the cell decomposition. My original cell decomposition involves infinitely many cells near a singular point.)

Since $g'(x)$ is either $>0$ or $<0$ over the entire interval $I$, $g$ has no critical point over $I$. For the inflection points, since $\frac{\partial f}{\partial x}$ and $\frac{\partial f}{\partial y}$ are both nonzero in $B(q) \setminus \{q\}$, equation (\ref{eqn:F2.1}) implies that 
\begin{equation*}
g'(x) =- \frac{f_x(x, g(x))}{f_y(x, g(x))}.
\end{equation*}
Substituting this into equation (\ref{eqn:F2.6}), we obtain an expression for $g''(x)$: 
\begin{equation*}
g''(x)= \frac{2f_x(x, g(x))f_y(x, g(x))f_{xy}(x, g(x)) - f_y(x, g(x))^2f_{xx}(x, g(x)) - f_x(x, g(x))^2f_{yy}(x, g(x))}{f_y(x, g(x))^3}.
\end{equation*}
It follows that if $g''(x)=0$, $(x, g(x))$ is contained in the following variety: 
\begin{equation*}
V(f) \cap V(2f_xf_yf_{xy} - f_y^2f_{xx} - f_x^2f_{yy}).
\end{equation*}
In the case that the polynomial $2f_xf_yf_{xy} - f_y^2f_{xx} - f_x^2f_{yy}$ is zero or is divisible by $f$, $g''$ is identically equal to zero thus having no strict inflection point. Otherwise, under the assumption that $f$ is irreducible, $f$ and $2f_xf_yf_{xy} - f_y^2f_{xx} - f_x^2f_{yy}$ are two nonzero polynomials in $\mathbb{R}[x,y]$ with no common factors, so the variety $V(f) \cap V(2f_xf_yf_{xy} - f_y^2f_{xx} - f_x^2f_{yy})$ is a finite set. Therefore, $g''(x)$ has at most finitely many strict inflection points over the open interval $I$.

\end{proof}

We finish the proof Theorem~\ref{thm:F2.1} as follows. According to Lemma~\ref{lem:F2.3}, a component $P$ in case (iii) is the graph of some real analytic function $g(x)$ over some open interval, then as before we can assign 0-cells to all strict inflection points and the two endpoints (one at the singular point $q$ and the other on the boundary of $B(q)$). Next, we can assign 1-cells to all line segments in between these finitely many adjacent 0-cells. Then, we repeat this same procedure for each component inside $B(q) \setminus \{q\}$, obtaining a finite cell decomposition for $X \cap B(q)$. (If there is no line component inside $B(q) \setminus \{q\}$, $X \cap B(q)$ is a single point at $q$. Assign a 0-cell at point $q$.)

Finally, since $X$ (by the compactness) can be covered by finitely many open sets in the form of either $I \times J$ around a non-singular point $p$ or $B(q)$ centered at a singular point $q$ in $S$, $X$ has a finite cell decomposition. 
It is easy to check that each cell in $\mathcal{A}$ is a semi-algebraic set. Indeed, any 1-cell is a continuous open line segment on the variety $V(f)$, and so is the intersection of an open rectangle with $V(f)$, which is semi-algebraic. 
If $\gamma$ is a shortest-length piecewise $C^2$-curve between two points in $X$, the intersection of $\gamma$ with each cell in $\mathcal{A}$ is either empty or consists of a single component that is either a singleton or a geodesic line segment. 
\end{proof}

Theorem~\ref{thm:F2.1} assumes that $X$ satisfies the compactness and connectedness properties, the next corollary shows that these two conditions are actually redundant. 

\begin{corollary}\label{cor:F2.2}
Suppose $f(x,y)$ is an irreducible polynomial function and $X = V(f)$ is the affine algebraic variety determined by the zeros of $f$. Then $X$ has a cell decomposition $\mathcal{A}$ with the desired finiteness property. 
\end{corollary}

\begin{proof}
It suffices to prove for the case when $X$ is connected, but not necessarily compact. In general, if each connected component of $X$ has a cell decomposition with the finiteness property, so does $X$. From now on, let us assume that $X$ is connected. 

There exists a large positive integer $N$ such that the open ball $B(0, N)$ centered at 0 with radius $N$ contains all points of $X$ for which either $\frac{\partial f}{\partial x} = 0$, or $\frac{\partial f}{\partial y} = 0$. Such an $N$ exists, because of (\ref{eqn:F2.2}) and (\ref{eqn:F2.3}). For the part of $X$ contained inside the closed ball $\bar{B}(0, N)$, it can again be covered by finitely many open sets in the form of either $I \times J$ around a non-singular point $p$ or $B(q)$ for a singular point $q$. It guarantees the existence of a cell decomposition with the finiteness property using the same proof as Theorem~\ref{thm:F2.1}.

Next, for each $n \geq N$, consider the closed annulus $\bar{A}(n, n+1)$ centered at 0 with inner and outer radii being $n$ and $n+1$, respectively. Then the intersection of $X$ with $\bar{A}(n, n+1)$ can be covered by finitely many open sets in the form of only $I \times J$, which also guarantees the existence of a cell decomposition with the finiteness property.

Lastly, we combine the cell decomposition for $X \cap \bar{B}(0, N)$ with those for $X \cap \bar{A}(n, n+1)$, where $n \geq N$, thus yielding a cell decomposition $\mathcal{A}$ with the desired finiteness property. 
\end{proof}

Now we are ready for the following general theorem concerning an arbitrary affine algebraic set in $\mathbb{R}^2$. 

\begin{theorem}\label{thm:F2.2}
Suppose $X$ is any arbitrary affine algebraic set in the plane, then $X$ has a cell decomposition $\mathcal{A}$ with the finiteness property. Furthermore, $(X, \mathcal{A})$ is a CW complex. 
\end{theorem}

\begin{proof}
Since $\mathbb{R}$ is Noetherian, Hilbert Basis Theorem implies that $\mathbb{R}[x,y]$ is also Noetherian. Then $X = V(f_1, \ldots, f_k)$ for finitely many polynomials $f_i \in \mathbb{R}[x,y]$, where $1~\leq~i~\leq~k$. 

For each $f_i$, we can write it as a product of finitely many irreducible polynomials, say $f_i^1, \ldots, f_i^{\, m_i}$. Since $V(f, g) = V(f) \cap V(g)$ and $V(f\cdot g) = V(f) \cup V(g)$, $X$ can be rewritten as follows:
\begin{equation*}
X= \bigcap_{i = 1}^k \ V(f_i^1) \cup \ldots \cup V(f_i^{\, m_i}).
\end{equation*}
Distributing the intersections over the unions, it turns out that $X$ is a finite union of sets in the following form:
\begin{equation*}
V(f_1^{j_1}) \cap V(f_2^{j_2}) \cap \ldots \cap V(f_k^{j_k}), \text{ where } 1 \leq j_i \leq m_i,  \text{ for each } i = 1, \ldots, k. 
\end{equation*}

If $k \geq 2$, the above expression consists of at most finitely many points because of the algebraic geometry theorem that we've utilized earlier \cite{P}. Therefore $X$ is either an empty set or consists of finitely many points in $\mathbb{R}^2$, so the theorem follows trivially. 

Next suppose that $k = 1$, then $X = V(f_1^1) \cup \ldots \cup V(f_1^{\, m_1})$. If $m_1=1$, we are done. If not, for each $j = 1, \ldots, m_1$, $V(f_1^j)$ has a cell decomposition $\mathcal{A}_j$ with the finiteness property based on Corollary~\ref{cor:F2.2}. Consider the following set $T$:
\begin{equation*}
T = \{x \in X: x \in V(f_1^{j}) \cap V(f_1^{j'}), \text{ where } 1 \leq j \neq j' \leq m_1\}.
\end{equation*}
Then $T$ is finite. Adjust $\mathcal{A}_j$ for each $j$ by adding a 0-cell for each point of $V(f_1^j)$ that lies in $T$, and then including extra 1-cells if necessary. Call the new cell decomposition $\mathcal{A}'_j$. 

Let $\mathcal{A} = \mathcal{A}'_1 \cup \ldots \mathcal{A}'_{m_1}$. If $\gamma$ is a shortest-length piecewise $C^2$-curve between two points in $X$, then 
\begin{equation*}
\gamma = ( \ \gamma \ \cap V(f_1^1) \ ) \  \bigcup \  \ldots \  \bigcup \  ( \ \gamma \ \cap V(f_1^{\, m_1}) \ ).
\end{equation*}
By the compactness of $\gamma$, for each $j = 1, \ldots, m_1$, $\gamma \ \cap V(f_1^j)$ consists of finitely many components, each of which is either a singleton or a shortest-length piecewise $C^2$-path between its two endpoints. Therefore, there are at most finitely many components in the intersection of $\gamma$ with each 1-cell in $\mathcal{A}$, each of which is either a singleton or a geodesic line segment. Moreover, if there are more than one component, then at least one endpoint of one of the components must lie in $T$, thus contradicting the fact that none of the 1-cells in $\mathcal{A}$ contains a point in $T$. As a conclusion, $\gamma$ intersects each 1-cell in $\mathcal{A}$ at most once, and the intersection must be a geodesic line segment. 

Since $X$ is a Hausdorff space and $\mathcal{A}$ is a cell decomposition, $(X, \mathcal{A})$ is a cell complex. Moreover, $\mathcal{A}$ is locally finite based on the construction in Corollary~\ref{cor:F2.2}. The two additional conditions (C) and (W) for a CW complex are automatic \cite{L2}. Therefore $(X, \mathcal{A})$ is a CW complex. 
\end{proof}

\begin{remark}
In the proof of Theorem~\ref{thm:F2.2}, we may take $\mathcal{A}$ to be the union of $\mathcal{A}_1, \ldots, \mathcal{A}_{m_1}$, directly. Then the intersection of $\gamma$ with every 1-cell in this cell decomposition may consist of more than one component, each of which is either a singleton or a geodesic line segment. As a consequence, this cell decomposition also works for the purpose of proving the theorem. However, we have chosen to adjust $\mathcal{A}_j$, for each $j$, in order to obtain a nicer cell decomposition, as illustrated in the proof of Theorem~\ref{thm:F2.2}. 
\end{remark}

\begin{corollary}\label{cor:F2.4}
Suppose $X$ is a finite union of arbitrary affine algebraic sets in the plane, then $X$ has a cell decomposition $\mathcal{A}$ with the finiteness property. Furthermore, $(X, \mathcal{A})$ is a CW complex. 
\end{corollary}

\begin{proof}
Since $ V(f) \cap V(g) = V(f^2+g^2)$ (communicated to me through Dr. Hardt), and $V(f) \cup V(g) = V(f \cdot g)$, there exists $f$ such that $X = V(f)$, which is thus an affine algebraic set. Applying Theorem~\ref{thm:F2.2} gives the desired result. 
\end{proof}

\section{The stratification of an open semi-algebraic set in $\mathbb{R}^2$}
Suppose $g_1(x, y), \ldots, g_m(x, y)$ are nonzero polynomial functions in two real variables, and assume they are irreducible and distinct. Define the open semi-algebraic set $Y$ as below:
\begin{equation} \label{eqn:F3.1}
Y = \{ g_1(x, y) > 0, \ldots, g_m(x, y) > 0 \}, \text{ where } m \geq 1. 
\end{equation}
Given such a $Y$, consider the following affine algebraic set $X$:
\begin{equation}\label{eqn:F3.2}
 X = V(g_1) \cup \ldots \cup V(g_m),
\end{equation}
which is closed. Then $\mathbb{R}^2 \setminus X$ is a disjoint union of connected open planar regions, in each of which the value of $g_j$ is either entirely greater than 0 or less than 0, due to the continuity of $g_j$ for each $j = 1, 2, \ldots, m$, and the connectedness. Therefore, $Y$ consists of some of (possibly none) these connected open planar regions. It suffices to come up with a proper stratification for each such individual open planar region, then a desired stratification for $Y$ is thus obtained by taking the union.

Since $X$ is also equivalent to $V(g_1 \cdots g_m)$, $X$ has a cellular stratification with the finiteness property based on Theorem~\ref{thm:F2.2}. Let's start with make an elementary observation.

\begin{lemma}\label{lem:F3.1}
Let $Y$ be given as in (\ref{eqn:F3.1}), and let $P$ be a nonempty connected component of $Y$. If the boundary of $P$ is nonempty, then it is contained in the affine algebraic set $X$ as defined in (\ref{eqn:F3.2}). 
\end{lemma}

\begin{proof}
Suppose $a$ is a boundary point of $P$, then $a$ is not inside $P$ otherwise it is an interior point. There exists a sequence $\{a_n\}_1^{\infty}$ of points in $P$ converging to $a$, which, by the continuity of $g_j$, implies that $g_j(a_n) \rightarrow g(a)$ for each $j = 1, 2, \ldots, m$. Thus $g_j(a) \geq 0$ for each $j$. 

If $g_j(a)=0$ for at least one $j$, then we are done. If not, $a$ is contained in one of the connected components in $Y$ other than $P$, making $a$ an exterior point of $P$.  
\end{proof}

From Lemma~\ref{lem:F3.1}, if $P$ is bounded, its boundary is nonempty (using the fact that $\mathbb{R}^2$ is connected and unbounded), and thus is contained in $X$. We want to employ the stratification of $X$ to get a cell decomposition for $P$. One such strategy is to divide $P$ using vertical strips whose endpoints are determined by the 0-cells on the boundary of $P$. 

\begin{figure}[ht]
\includegraphics[width=8cm]{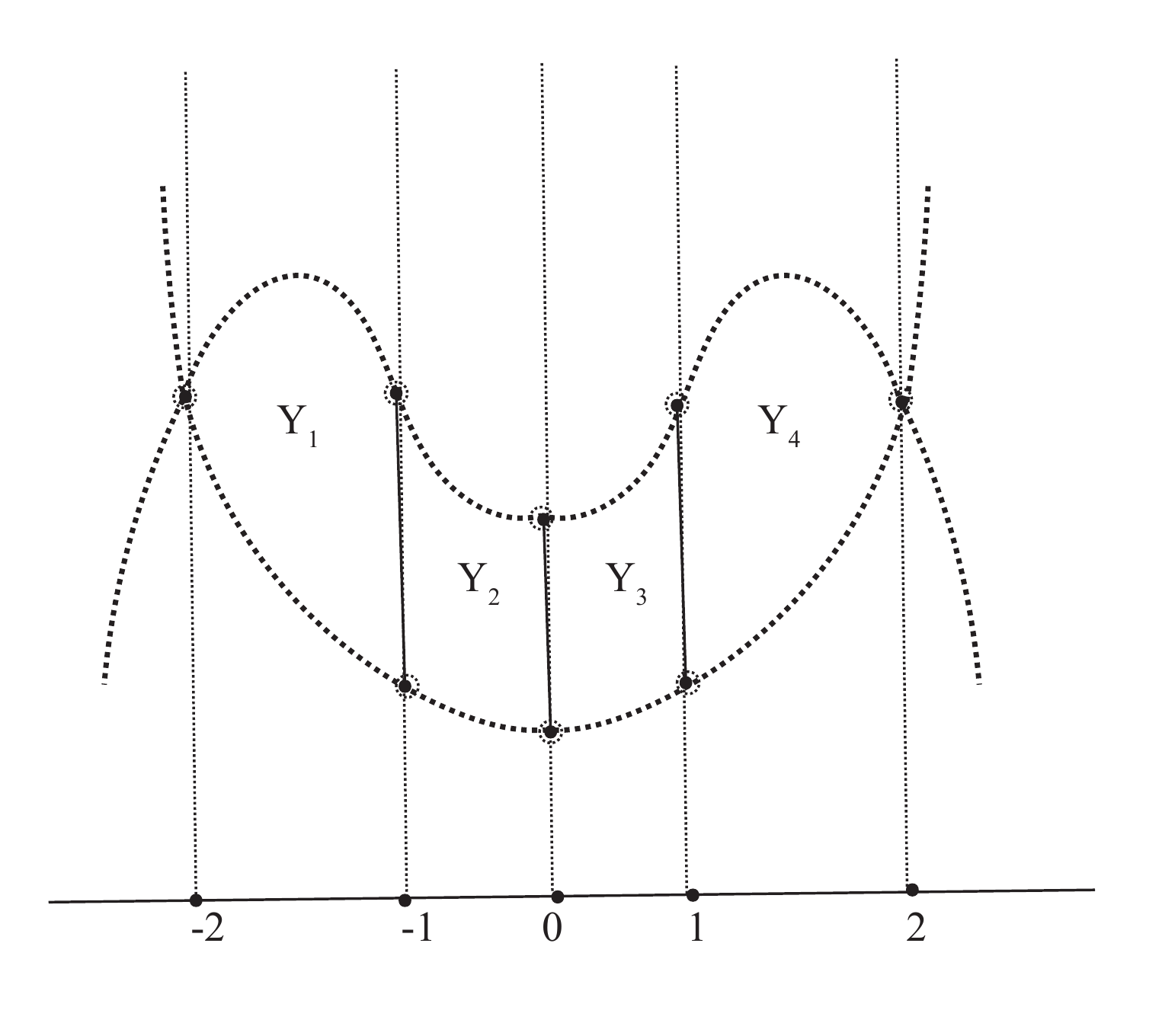}
\caption{Divide $Y$ through vertical strips determined by the 0-cells on the boundary of $Y$.}
\label{fig:F3.1}
\end{figure}

\begin{example}
let $g_1 = -y+(x^2-1)(x^2 -4)$ and $g_2=y-2(x-2)(x+2)$. Then $Y$ is the connected region bounded between two graphs as shown in Figure~\ref{fig:F3.1}. The 0-cells on the upper boundary of $Y$ determined by $g_1$ are at points $x = -2, -1, 0, 1, 2$; and the 0-cells on the lower boundary determined by $g_2$ are at points $x = -2, 0, 2$. Projecting these 0-cells upon the $x$-axis partitions $[-2,2]$ into five subintervals, each of which gives rise to a vertical strip. It follows that we can divide $Y$ into four sets: $Y_1, \ldots, Y_4$, each of which has a top lying on the graph of $g_1$, a bottom that is on the graph of $g_2$, and two sides being either a vertical line interval, or empty. We know how to construct a cell decomposition for each of the $Y_i$ from our previous discussions, thus obtaining a cell decomposition for $Y$. Indeed, for each $i=1, \ldots, 4$, $Y_i$ is a union of a region of type II and one or two open vertical line intervals. 
\end{example}

In general , we can apply a similar idea to the open semi-algebraic set $Y$. 
\begin{proposition}\label{prop:F3.2}
Let $Y$ be an open semi-algebraic set given as before, and $P$ be a connected component of $Y$. Suppose that $P$ is bounded, then $P$ has a cell decomposition $\mathcal{A}$ with the finiteness property. Furthermore, ($P$, $\mathcal{A}$) is a CW complex.
\end{proposition}


\begin{proof}
By hypothesis, $g_1, \ldots, g_m$ are irreducible and distinct, therefore the point of intersection of $V(g_i)$ and $V(g_j)$, for $i \neq j$, is at most a finite set. Include these points as 0-cells in the cell decomposition for each $V(g_i)$, $i = 1, \ldots, m$. (If a new 0-cell is within a 1-cell, divide the 1-cell into two new 1-cells.)
Indeed, since each cell decomposition for $V(g_i)$ exists by Theorem~\ref{thm:F2.2}, a cell decomposition for $X$ also exists by combining them. Call it $\mathcal{A}$, then ($X$, $\mathcal{A}$) is a CW complex because $\mathcal{A}$ is locally finite.

We claim that the cell decomposition $\mathcal{A}$ satisfies the property that every 0- or 1-cell in $\mathcal{A}$ is either entirely contained in the boundary of $P$ or entirely not. This is obviously true for all the 0-cells. For the 1-cells, the proof is as follows. 
Let $e_1$ be a 1-cell in $\mathcal{A}$ which has a nonempty intersection with the boundary of $P$. Without loss of generality, we may assume that $e_1$ belongs to the cell decomposition of $V(g_1)$. Consider the set 
\begin{equation*}
S = \{(x,y) \in e_1: \text{ $(x,y)$ is on the boundary of $P$ }\}, 
\end{equation*}
then $S$ is closed and nonempty. It suffices to show that $S$ is also open so that $S$ is equal to $e_1$ by the connectedness of $e_1$. 

Given $(x_0, y_0) \in S$, there exists an open neighborhood $I \times J$ around $(x_0, y_0)$ such that $\{(x, y) \in I \times J: g_1(x, y)=0\}$ is the graph of a real analytic function $h$ over the $x$-axis (or the $y$-axis). We may choose $I \times J$ to be so small that it doesn't intersect $V(g_2), \ldots, V(g_m$) due to the fact that $(x_0, y_0)$ is not an intersection point and so is at a positive distance from each of the closed sets $V(g_j)$, where $2 \leq j \leq m$. It follows that $g_j(x, y)$ is entirely $> 0$ for $(x, y) \in I \times J$ for each $2 \leq j \leq m$, since $g_j(x_0, y_0) > 0$ for all $j$. 
On the other hand, the graph of $h$ inside the open rectangle $I \times J$ divides it into two connected components, namely 
\begin{equation*}
A=\{(x, y) \in I \times J: y < h(x)\}  \text{ and  } B=\{(x, y) \in I \times J: y > h(x)\}.
\end{equation*}
It is easy to check that $g_1>0$ throughout at least one of the two components, pick the one that has a nonempty intersection with $P$, say $A$. Then, $A$ is a subset of $P$. As a result, every point in the set $\{(x, y) \in I \times J: g_1(x, y)=0\}$ is a boundary point of $P$. Therefore $S$ is also an open subset of $e_1$. Thus $e_1$ belongs to the boundary of $P$ and so does its closure. 
Since the boundary of $P$ is a subset of $X$ according to Lemma~\ref{lem:F3.1}, it follows that the boundary of $P$ is a finite subcomplex of $X$, because it is compact \cite{L2}.

Projecting the closure of $P$ onto the $x$-axis, the image is a finite closed interval, say $[a, b]$, where $-\infty < a < b < \infty$. Furthermore, projecting the 0-cells on the boundary of $P$ divides $[a, b]$ into finitely many intervals, say $a=a_0 < a_1 < \ldots < a_n = b$. 
For each $i = 1, \ldots, n-1$, we claim that the vertical line at $a_i$ intersects the boundary of $P$ at finitely many points. This is because $g_j(a_i, y) = 0$ has at most finitely many solutions for all $j = 1, \ldots, m$, unless $g_j(a_i, y) \equiv 0$. If $g_j(a_i, y) \equiv 0$ for some $j$, then the vertical line at $a_i$ divides the plane into two halves, so $P$ lies inside only one of the two halves. Therefore, $a \geq a_i$ or $b \leq a_i$, resulting in a contradiction. On the other hand, at $x=a$ or $x=b$, the intersection of the vertical line with the boundary of $P$ is a finite disjoint union of closed vertical line intervals (of finite lengths) and isolated points. Indeed, the intersection is compact, so it consists of only finitely many connected components, each of which is a connected subset of a real line. By the connectedness of the real line $\mathbb{R}$, if a connected component has at least two points, it is an interval which is also closed and bounded in our case; if a connected component has only one point, then it is isolated from the others with respect to the subspace topology induced from $\mathbb{R}$. 

For each $i = 1, \ldots, n-1$, add the finitely many points of intersection of the vertical line at $a_i$ and the boundary of $P$ as 0-cells in the cell decomposition for the boundary of $P$. For $i = 0$, or $n$, we add these new 0-cells: the finitely many isolated points and the {\it endpoints} of the vertical line segments in the intersection of the vertical line at $a$ or $b$ with the boundary of $P$. It follows that with respect to this new cell decomposition, each 1-cell on the boundary of $P$ lies directly over one and only one open interval $(a_i, a_{i+1})$ for some $i = 0, \ldots, n-1$, except for the possible vertical 1-cells at the two endpoints $a$, $b$ and for those 1-cells that are graphs over the y-axis. 

We want to modify these 1-cells which are graphs over the y-axis so that they become graphs over the x-axis as well. This can be done by dividing these 1-cells further. Let $e_2$ be one of such 1-cells, and without loss of generality, we may assume that $e_2$ is carried by $V(g_1)$. Based on our construction in Theorem~\ref{thm:F2.1}, there exists an open neighborhood $I_2 \times J_2$ of $e_2$ such that 
\begin{equation*}
\bar{e}_2 = \{(x, y) \in \bar{I}_2 \times \bar{J}_2: g_1(x, y)=0\} = \{(h_2(y), y): y \in \bar{J}_2\}, 
\end{equation*}
for some real analytic function $h_2$. If $h_2$ is linear, it's either constant, in which case we get a vertical line segment, or has a nonzero slope, in which case the graph of $h$ as a function of $y$ is also a graph over the $x$-axis. Suppose $h_2$ is nonlinear, we can insert its local {\it maximum} points as 0-cells, which are finitely many according to Corollary~\ref{cor:F2.3}. It follows that every new 1-cell in $e_2$ is either strictly increasing or decreasing over the y-axis, thus becoming a graph over the x-axis as well. Furthermore, each such graph when viewed as over the $x$-axis is also real analytic, because given $(x_1, y_1) \in e_2$ with $h'(y_1) \neq 0$, 
\begin{equation*}
\frac{\partial g_1}{\partial x}(x_1, y_1) \neq 0, 
\end{equation*}
by hypothesis, and moreover 
\begin{equation*}
\frac{\partial g_1}{\partial x}(x_1, y_1) \cdot h'(y_1) + \frac{\partial g_1}{\partial y}(x_1, y_1) =0 \Longrightarrow \frac{\partial g_1}{\partial y}(x_1, y_1) \neq 0.
\end{equation*}

Repeat the previous process for these new 0-cells, that is, first project them onto the $x$-axis, then include as 0-cells for the points of intersection of the vertical lines with the boundary of $P$. In the end, we obtain a cell decomposition in which each 1-cell is either an open vertical interval or a real analytic function over $(a_i, a_{i+1})$ for some $i = 0, \ldots, n-1$ (using the same notation for the new partition of $[a, b]$ as before). Furthermore, there are at least two separate non-vertical 1-cells lying over $(a_i, a_{i+1})$ for all $i$, otherwise $P$ is disconnected or unbounded. A picture for the intersection of the closure of $P$ with the vertical strip $[a_i, a_{i+1}] \times \mathbb{R}$ is shown in Figure~\ref{fig:F3.2}. 
We note that the 1-cells lying over a common interval might share common endpoints in their closures. 

\begin{figure}[ht]
\includegraphics[width=12cm]{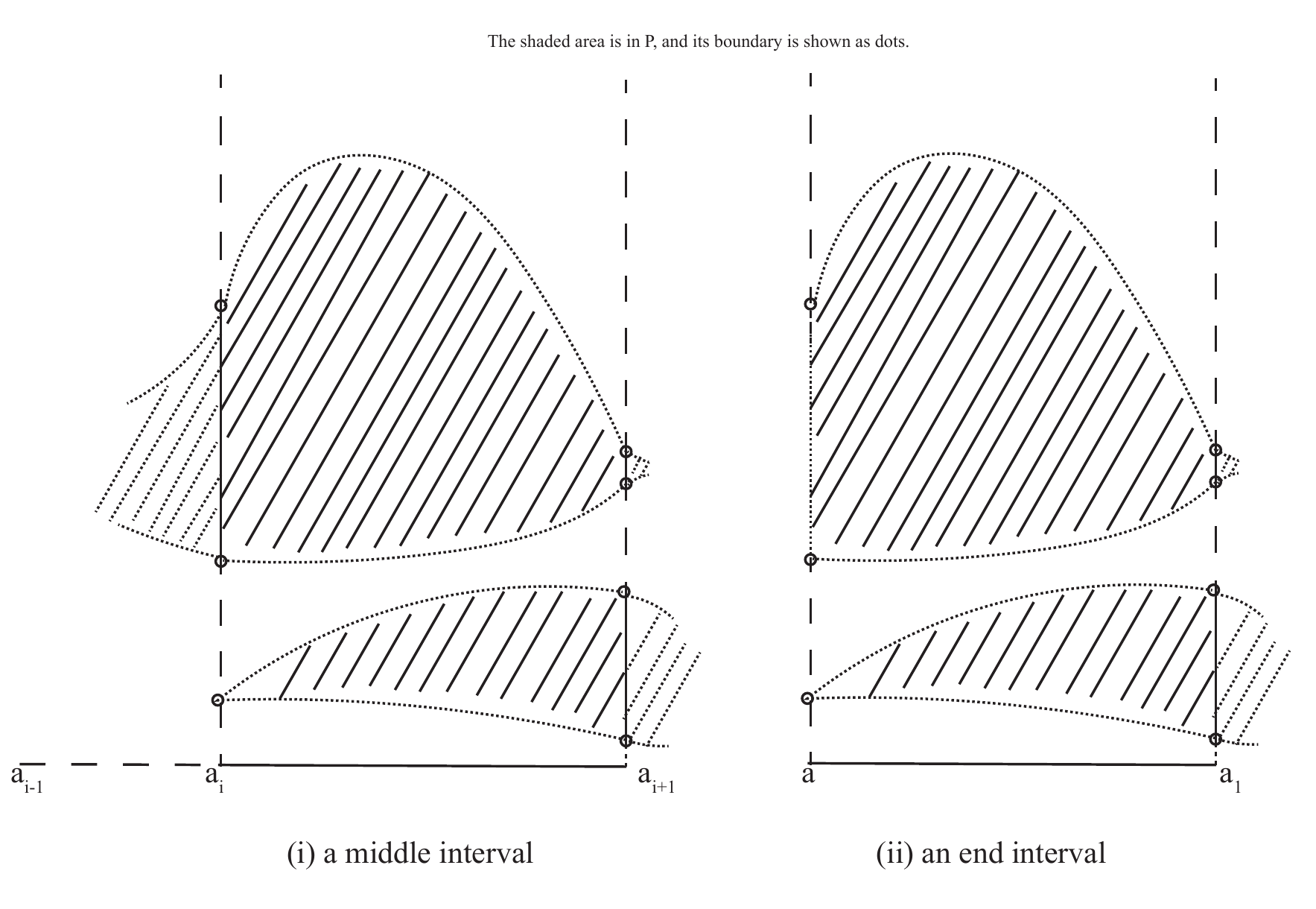}
\caption{The part of $P$ and its boundary lying over an interval $[a_i, a_{i+i}]$.}
\label{fig:F3.2}
\end{figure}

For each $(a_i, a_{i+1})$, since only finitely many 1-cells spread out over it, the open vertical strip $(a_i, a_{i+1}) \times \mathbb{R}$ subtracting these 1-cells consist of finitely many connected open regions bounded by at least one of these 1-cells on the top or bottom, and by a vertical line interval or a point on the two sides (see Figure~\ref{fig:F3.3}). We call such a connected region {\bf a basic (open) region}. (An exception will be discussed soon.)
Then the intersection of $P$ and the open vertical strip is a finite union of these bounded basic regions which intersect $P$ nontrivially. 
It follows that $P$ can be partitioned into finitely many basic regions together with some of the open vertical sides for these regions. Two examples have already been shown in Figure~\ref{fig:F3.1} and Figure~\ref{fig:F3.2}. 
\begin{figure}[ht]
\includegraphics[width=10cm]{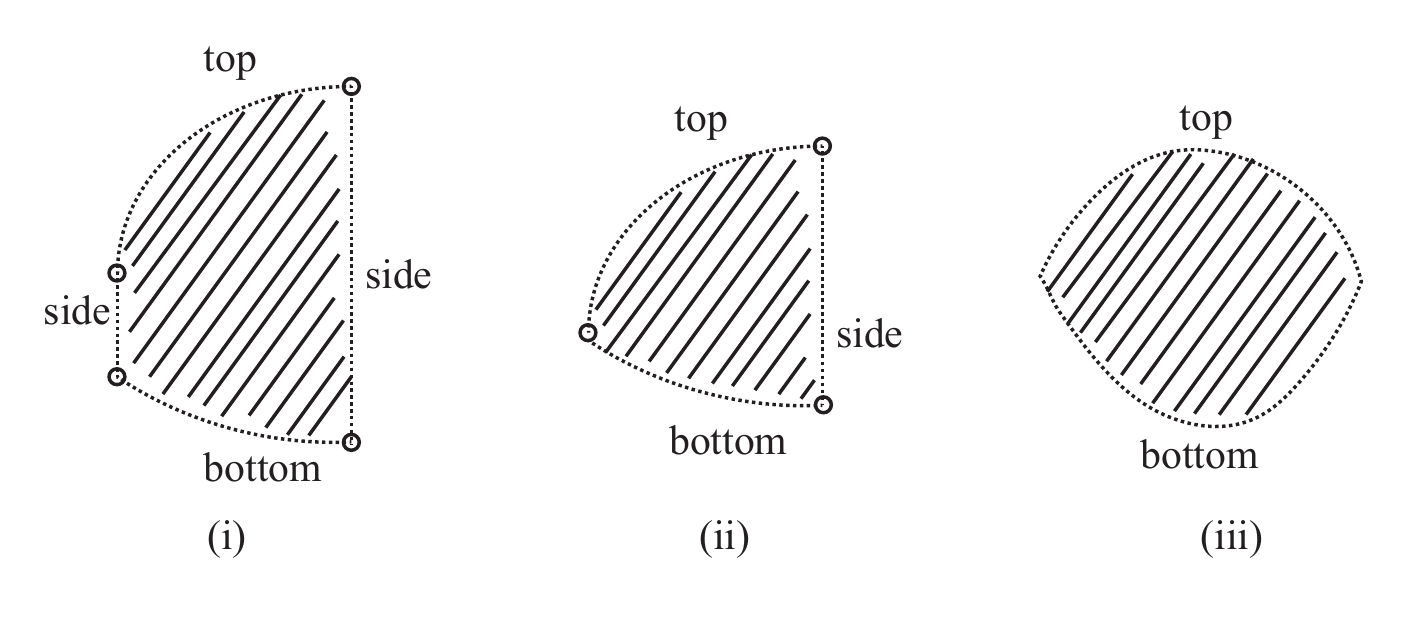}
\caption{A basic (open) region.}
\label{fig:F3.3}
\end{figure}

Since a basic region is a region of type II, a cell decomposition with the finiteness property exists \cite{Y}. 
More precisely, each basic region is a finite intersection of open polynomial half planes (that is, an open region below the graph of a polynomial function up to rotation, which is analogous to an open half plane); each polynomial half plane has a cell decomposition with a sequence of 1-cells (in the same shape of the graph) converging to the boundary; then overlaying the cell decomposition for each open polynomial half plane leads to a cell decomposition for the basic region. 

For any vertical open interval, whenever it is included as a side of a basic region in $P$, we
don't do anything with that side. That is to say, we exclude the cell decomposition of an open polynomial half plane corresponding to that particular side when performing the overlapping. 
One can check that such a cell decomposition is actually locally finite. Since $P$ is also Hausdorff, $(P, \mathcal{A})$ is a CW complex. 

There is one exception that we need to discuss carefully. The cell decomposition for $X$ in (\ref{eqn:F3.2}) may consist of isolated points. (For example, $x^2+y^2=0$). Therefore, on each side of a basic region that is a vertical line interval, there might have an isolated point removed as shown in Figure~\ref{fig:F3.9}. 
\begin{figure}[ht]
\includegraphics[width=5cm]{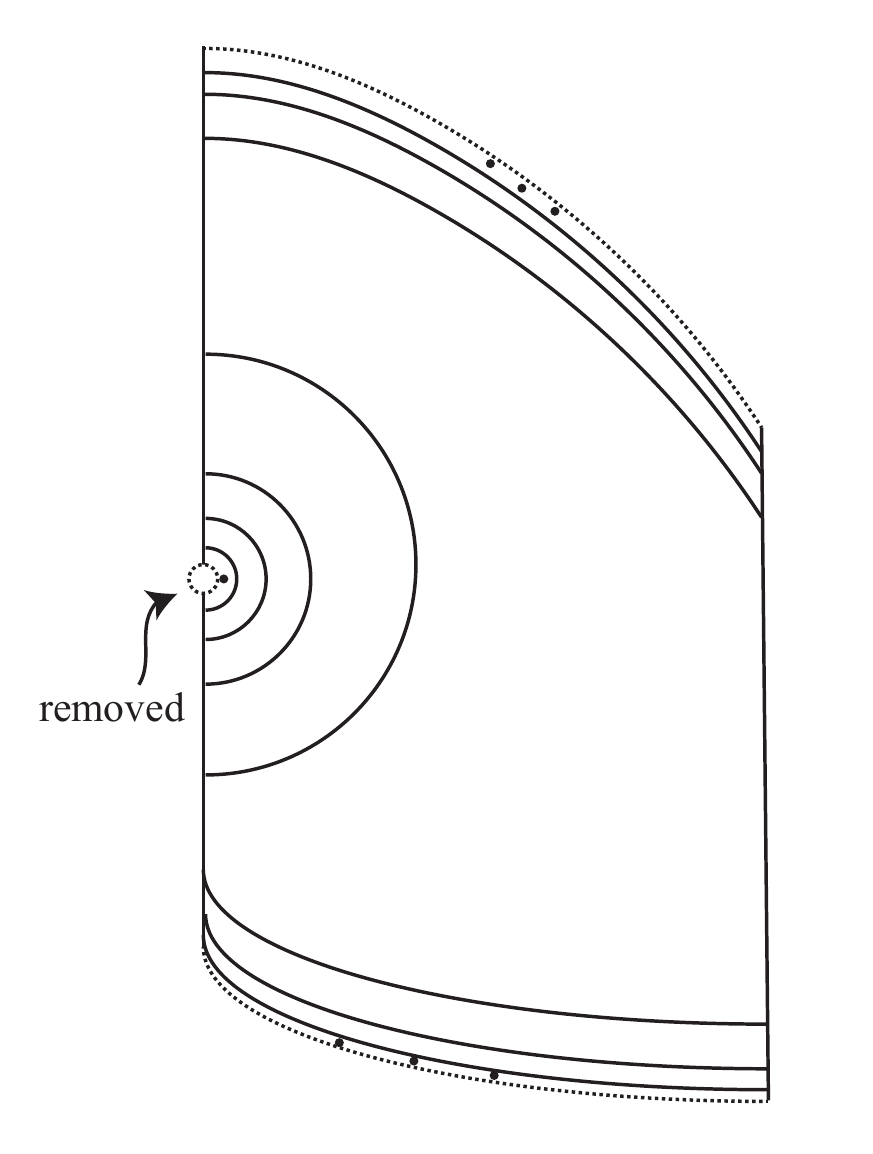}
\caption{The cell decomposition of a basic region with an isolated point removed on one of its two sides.}
\label{fig:F3.9}
\end{figure}
In such a situation, let's use semicircles to divide the basic region even more so that the local finiteness property can be achieved. It suffices to look at the following example as an illustration.

\begin{example}
Let $C = [-1,1]\times[0,1] \setminus \{(0,0)\}$ be a cube in the plane with the origin being removed, then we want to come up with a cell decomposition of $C$ that is also a CW complex. Use a sequence of semicircles centered at $(0,0)$ with radii $1/n$, $n \geq 1$, to divide $C$. It follows that such a cell decomposition is locally finite (see Figure~\ref{fig:F3.9}).
\end{example}

Furthermore, we won't have a situation in which for the two 1-cells adjacent to the isolated point being removed, one is in $P$ and the other is not in $P$. See Figure~\ref{fig:F3.10}.
\begin{figure}[ht]
\includegraphics[width=4cm]{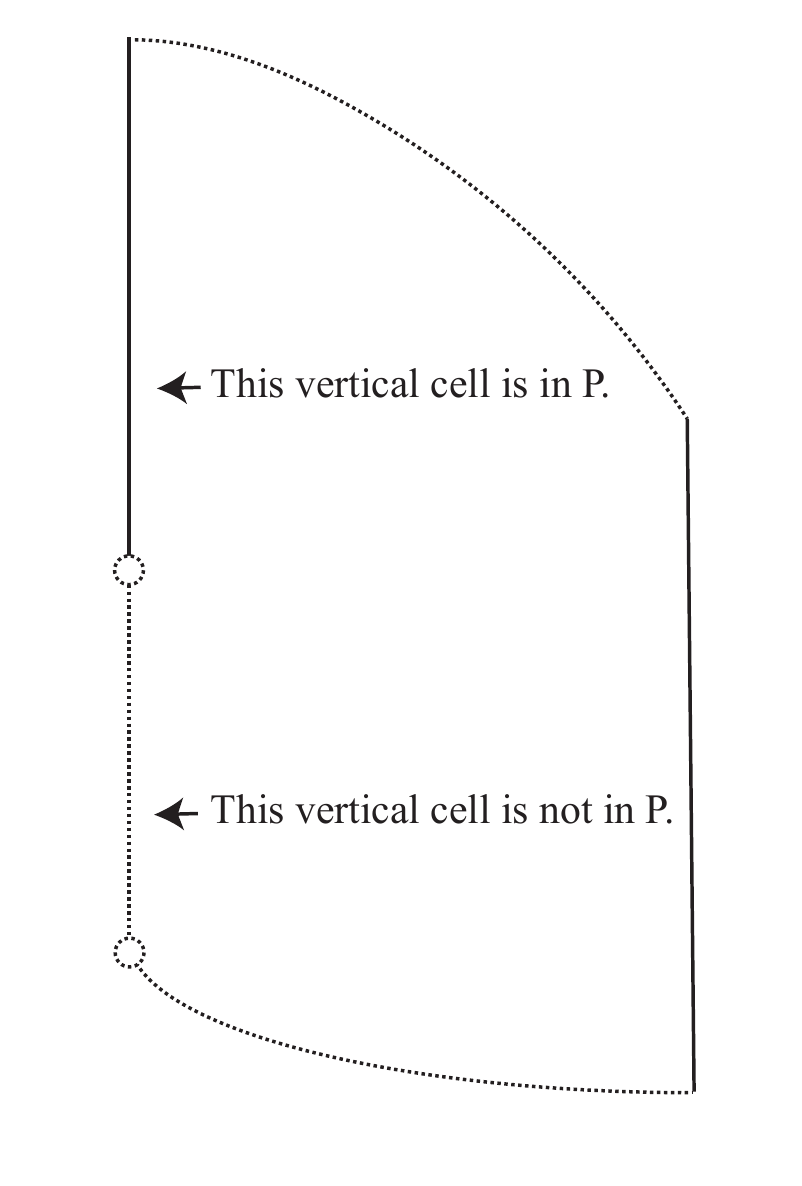}
\caption{It is impossible to have a basic region with an isolated point removed such that the two adjacent 1-cells are not simultaneously inside or outside $P$.}
\label{fig:F3.10}
\end{figure}
If one 1-cell is not in $P$, it must lie in $X$, because $g_1 \geq 0, \ldots, g_m \geq 0$ on it due to continuity. Then this 1-cell is over the endpoint $a$ or $b$, and so one of the $g_j = 0$ is a vertical line at $a$ or $b$. Therefore the other 1-cell is also not in $P$. On the other hand, suppose one of the two 1-cells is in $P$, then they must correspond to an $a_i \neq a, b$. Since $x = a_i$ intersects the boundary of $P$ at finitely many points, the other cell must also be in $P$. Thus the situation in Figure~\ref{fig:F3.10} will never happen. 

It is easy to check that each cell in $\mathcal{A}$ is a semi-algebraic set. Indeed, each 2-cell is bounded by finitely many 1-cells each of which is either a vertical/horizontal line interval, or inherits the same shape from one of the top and bottom graphs of a basic region, or is a semicircle. It follows that all these 1-cells can be determined by algebraic equations, therefore making the 2-cell a semi-algebraic set.  

Given a shortest-length curve $\gamma$ between two points in $P$ (if it exists), $\gamma$ intersects the closure of each 1-cell at most finitely many times. Indeed, the closure of each 1-cell is the graph of some real analytic function. If $\gamma$ intersects it infinitely many times, we can find an accumulation point in the intersection. Since $\gamma$ is locally a straight line ($P$ is an open set), the 1-cell must be also linear thus leading to a contradiction. Therefore, each 1-cell interacts $\gamma$ at most finitely many times, thus so does every 2-cell. As a result, the cell decomposition $\mathcal{A}$ satisfies the desired finiteness property. 
This finishes the proof of the proposition. 
\end{proof}

In the proposition, we assume that $P$ is a bounded component of $Y$. In fact, this hypothesis can be removed by dividing the plane into cubes and focusing the cell decomposition in each cube. This idea is formally stated as follows: 

\begin{proposition}\label{prop:F3.3}
Suppose $Y$ is an open semi-algebraic set defined as in (\ref{eqn:F3.1}), and let $P$ be a connected component of $Y$. Then there exists a cell decomposition $\mathcal{A}$ for $P$ such that $(P, \mathcal{A})$ is a CW complex satisfying the finiteness property. 
\end{proposition}

\begin{proof}
If $P$ is equal to the whole plane, then there is nothing to prove. Otherwise, the boundary of $P$ is nonempty. 
Let's pick a positive integer $R$ large enough so that the four sides of the cube $[-R, R] \times [-R, R]$ intersects the boundary of $P$ at mostly finitely many points. Indeed, if there were infinitely many intersection points, say on the side $x = R$, then there exists $1 \leq j \leq m$ such that $V(g_j) \cap \{x=R\}$ is an infinite set. Since both $g_j$ and $x-R$ are irreducible, they must have a common factor implying that $g_j = c(x-R)$ for some nonzero constant $c$. In this situation, enlarging $R$ fixes the problem. In fact, if $g_j =c( x-a)$, or $c(y-a)$ for some $a \in \mathbb{R}$, where $1 \leq j \leq m$, we require $R$ to be bigger than $a$. 

Suppose the cube $[-R, R] \times [-R, R]$ intersects the boundary of $P$ trivially, then the cube is either entirely contained in $P$ or entirely not, due to the connectedness property of the cube. In the case that it is completely inside $P$, a cell decomposition for it is shown in Figure~\ref{fig:F3.5}. This is quite obvious. The cells compose of the four vertices, the interiors of the four sides, and the interior of the cube. 

\begin{figure}[ht]
\includegraphics[width=3cm]{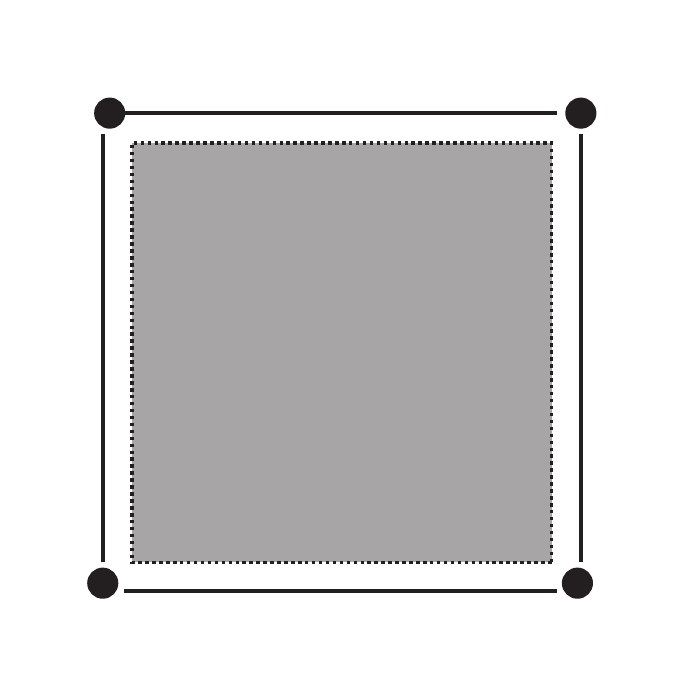}
\caption{A cell decomposition of a cube in $\mathbb{R}^2$.}
\label{fig:F3.5}
\end{figure}

Suppose the cube $[-R, R] \times [-R, R]$ intersects the boundary of $P$ nontrivially. The boundary of $P$ inside the cube contains finitely many 0-cells. Like in Proposition~\ref{prop:F3.2}, we add at most finitely many more 0-cells to ensure that every 1-cell is a graph of some real analytic function over the $x$-axis or an open vertical interval. Moreover, we want to also include the following points as 0-cells: the four vertices of the cube, and
the points of intersection of the boundary of $P$ and the {\it boundary} of the cube. Projecting these 0-cells onto the $x$-axis partitions $[-R, R]$ into finitely many subintervals, say $-R = a_0 < a_1 < \ldots < a_n=R$. 
Over each such subinterval $[a_i, a_{i+1}]$, the bounded vertical strip $(a_i, a_{i+1}) \times [-R, R]$ is again being subdivided into finitely many connected components, each of whose interior is either contained completely inside $P$ or not (see Figure~\ref{fig:F3.6}). Each such an interior region is a basic (open) region as we've seen earlier. Therefore, the intersection of $P$ with the closed cube $[-R, R] \times [-R, R]$ is a finite disjoint union of basic regions, and finitely many vertical or horizontal intervals which could be open, closed, or half-closed. It follows that a cell decomposition exists with the finiteness property. Indeed, if any vertical or horizontal open interval is included as a side of a basic region in $P$, we do nothing with that side as seen in Proposition~\ref{prop:F3.2}. Moreover, if an endpoint of a vertical or horizontal open interval is inside $P$, that side must be entirely in $P$ as well. Indeed, every vertical or horizontal open interval lies completely inside either $P$, or the boundary of $P$, or the exterior of $P$. Therefore, if a side does not lie in $P$, then it is inside the complement of $P$, which is a closed set thus having the two endpoints of the side as well. It follows that this observation guarantees that our cell decomposition satisfies the local finiteness property. 

\begin{figure}[ht]
\includegraphics[width=10cm]{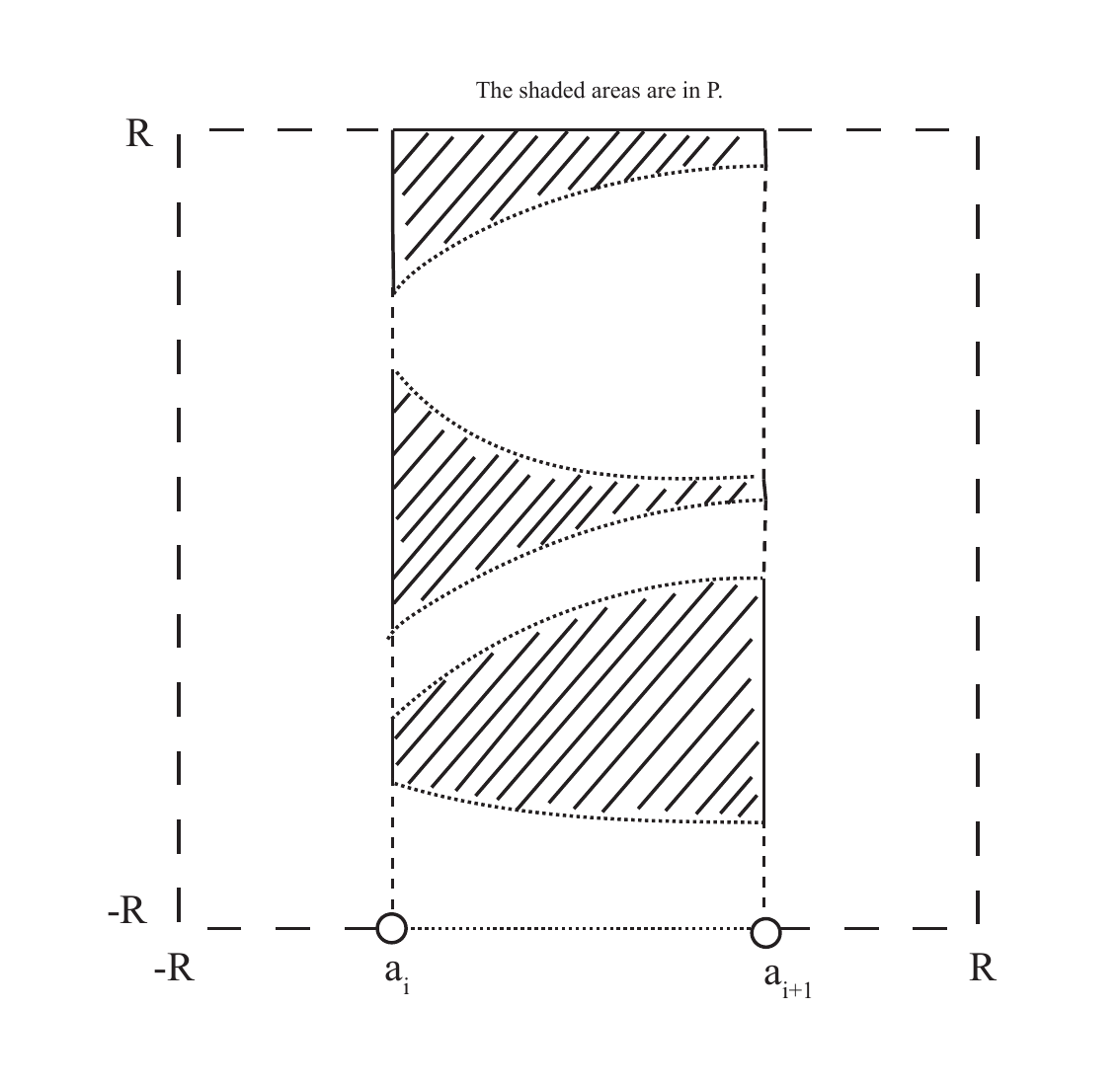}
\caption{The vertical strip $(a_i, a_{i+1}) \times [-R, R]$ and the 1-cells contained in it. }
\label{fig:F3.6}
\end{figure}

Since the plane $\mathbb{R}^2$ is a countable union of the cube $[-R, R] \times [-R, R]$ and other cubes in the form of 
\begin{equation*}
[nR, (n+1)R] \times [mR, (m+1)R],
\end{equation*}
where at least one of  $n, m \geq 1,$ or $\leq -2$ (See Figure~\ref{fig:F3.13}).
For each one of these cubes, its boundary intersects the boundary of $P$ at mostly finitely many times. Applying the similar argument as for $[-R, R] \times [-R, R]$ gives a cell decomposition for the intersection of $P$ with each of these cubes. For two adjacent cubes, their common edge could get 0-cells from both cell decompositions, and this is fine because there are only finitely many 0-cells in total. Consequently, we obtain a cell decomposition $\mathcal{A}$ for $P$ which fulfills the local finiteness property, thus making $(P, \mathcal{A})$ a CW complex. Furthermore, $\mathcal{A}$ also satisfies the finiteness property. 
\begin{figure}[ht]
\includegraphics[width=6cm]{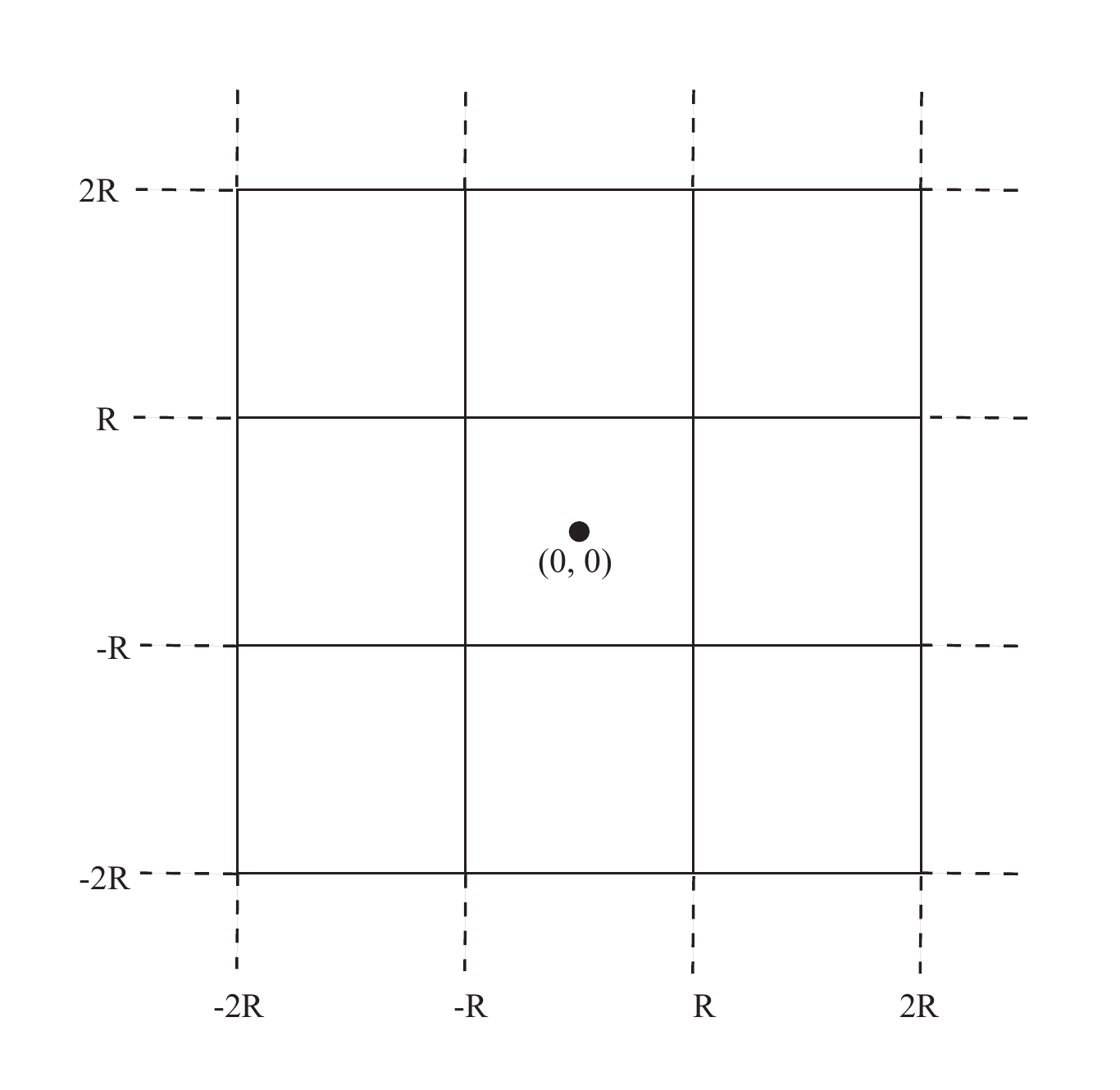}
\caption{When $P$ is unbounded, we may divide the plane into small cubes.}
\label{fig:F3.13}
\end{figure}
\end{proof}

\begin{theorem}\label{thm:F3.4}
Suppose $Y$ is an open semi-algebraic set defined as in (\ref{eqn:F3.1}), then there exists a cell decomposition $\mathcal{A}$ for $Y$ such that $(Y, \mathcal{A})$ is a CW complex satisfying the finiteness property.
\end{theorem}

\begin{proof}
Since this is true for each connected component of $Y$ according to the previous proposition, this is also true for $Y$. 
\end{proof}

Next we need to look at the general case of taking a finite union of sets in the form of $\{g_1(x, y) > 0, \ldots, g_m(x,y) > 0\}$, in which the $g_j$ are not necessarily irreducible. Since $g(x,y)\cdot h(x,y) > 0$ implies that either $g(x,y) >0, h(x,y)>0$ or $g(x,y) <0, h(x,y)<0$ (equivalently, $-g(x,y) > 0, -h(x,y) >0$). We may assume without loss of generality that the $g_j$ are indeed irreducible. 

\begin{theorem}\label{thm:F3.5}
Suppose $Y$ is a finite union of open semi-algebraic sets defined as in (\ref{eqn:F3.1}), more precisely, let $Y$ be
\begin{eqnarray*}
Y = & \{g_1(x, y) > 0, \ldots, g_m(x,y) > 0\} \cup \{g'_1(x, y) > 0, \ldots, g'_{m'}(x,y) > 0\} \cup \\  & \ldots  \cup \{g''_1(x, y) > 0, \ldots, g''_{m''}(x,y) > 0\},
\end{eqnarray*}
where the union is finite, and the $g_j$, $g'_j$, \ldots, $g''_j$ are all irreducible polynomials. 
Then there exists a cell decomposition $\mathcal{A}$ for $Y$ such that $(Y, \mathcal{A})$ is a CW complex satisfying the finiteness property.
\end{theorem}

\begin{proof}
The idea is similar as to the proof in Proposition~\ref{prop:F3.2}, however, there is a slight improvement as regard to the finiteness property. 

First, let's again definite $X$ to be the union of all boundaries as follows:
\begin{eqnarray*}
X  & =& \{g_1= 0\}, \ldots, \{g_m = 0\} \cup \{g'_1 = 0\}, \ldots, \{g'_{m'} = 0\} \cup \\  & &\ldots  \cup \{g''_1 = 0\}, \ldots, \{g''_{m''} = 0\},
\end{eqnarray*}
for which we can find a CW decomposition according to Corollary~\ref{cor:F2.3}. 

Second, look at each connected component $P$ that is in the union $Y$. Then the boundary $\partial P$ of $P$ is contained in $X$, and each 0- or 1- cell in $X$ is either entirely contained in $\partial P$ or not. Then chop $P$ up as before into basic regions (Propositions~\ref{prop:F3.2} and \ref{prop:F3.3}). Observe that if a 1-cell is not in $Y$, then its two endpoints are also not in $Y$, because the complement of $Y$ is closed. Thus the finiteness property has no problem for the 1-cells. However, this no longer holds for the 0-cells. We might have a 0-cell that is not in $Y$, but is adjacent to two 1-cells which are in $Y$ (For example, see Figure~\ref{fig:F3.7}).
\begin{figure}[ht]
\includegraphics[width=8cm]{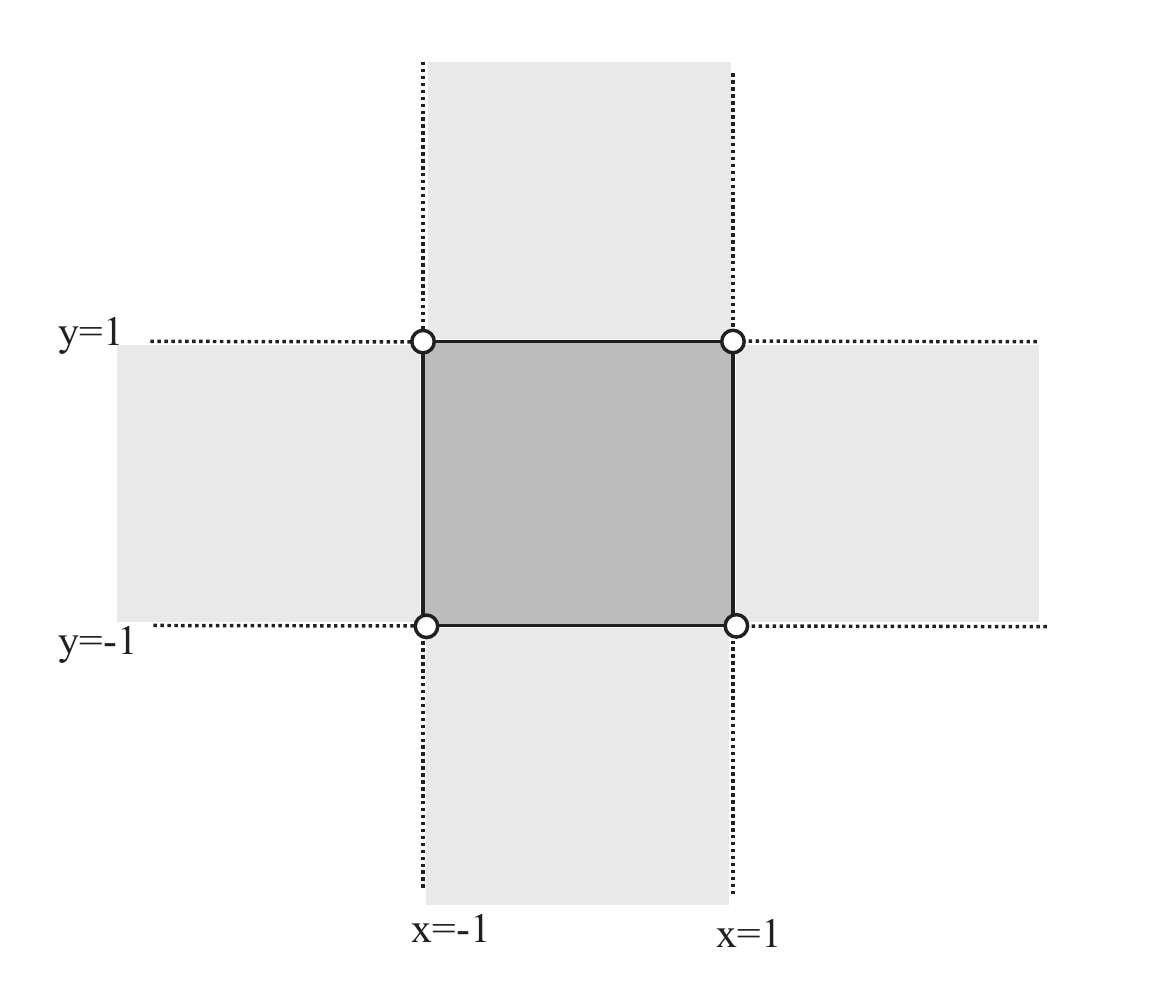}
\caption{When $Y = \{x>-1, x<1\} \cup \{y>-1, y<1\}$, look at the basic region in the middle, the four 0-cells at corners are not in $Y$, but they are adjacent to 1-cells that are in $Y$.}
\label{fig:F3.7}
\end{figure}
In order for the finiteness property to be satisfied near such a 0-cell, let's first look at the following example. 

\begin{example}
Suppose $C = [0,1]\times[0,1] \setminus \{(0,0)\}$ is the unit cube without one of its corners at $(0,0)$, then we want to look for a cell decomposition for $C$ that is also locally finite. In order to do so, let us first divide $C$ into four small cubes. Next, for the lower left cube which contains $(0,0)$, let's divide it further into another four small cubes. Then, pick the lower left cube which contains $(0,0)$, and divide it again. Continue this process infinitely many times. This actually results in a cell decomposition of $C$ that is also local finite everywhere in $C$. A demonstration is shown in Figure~\ref{fig:F3.8}. 
\begin{figure}[ht]
\includegraphics[width=6cm]{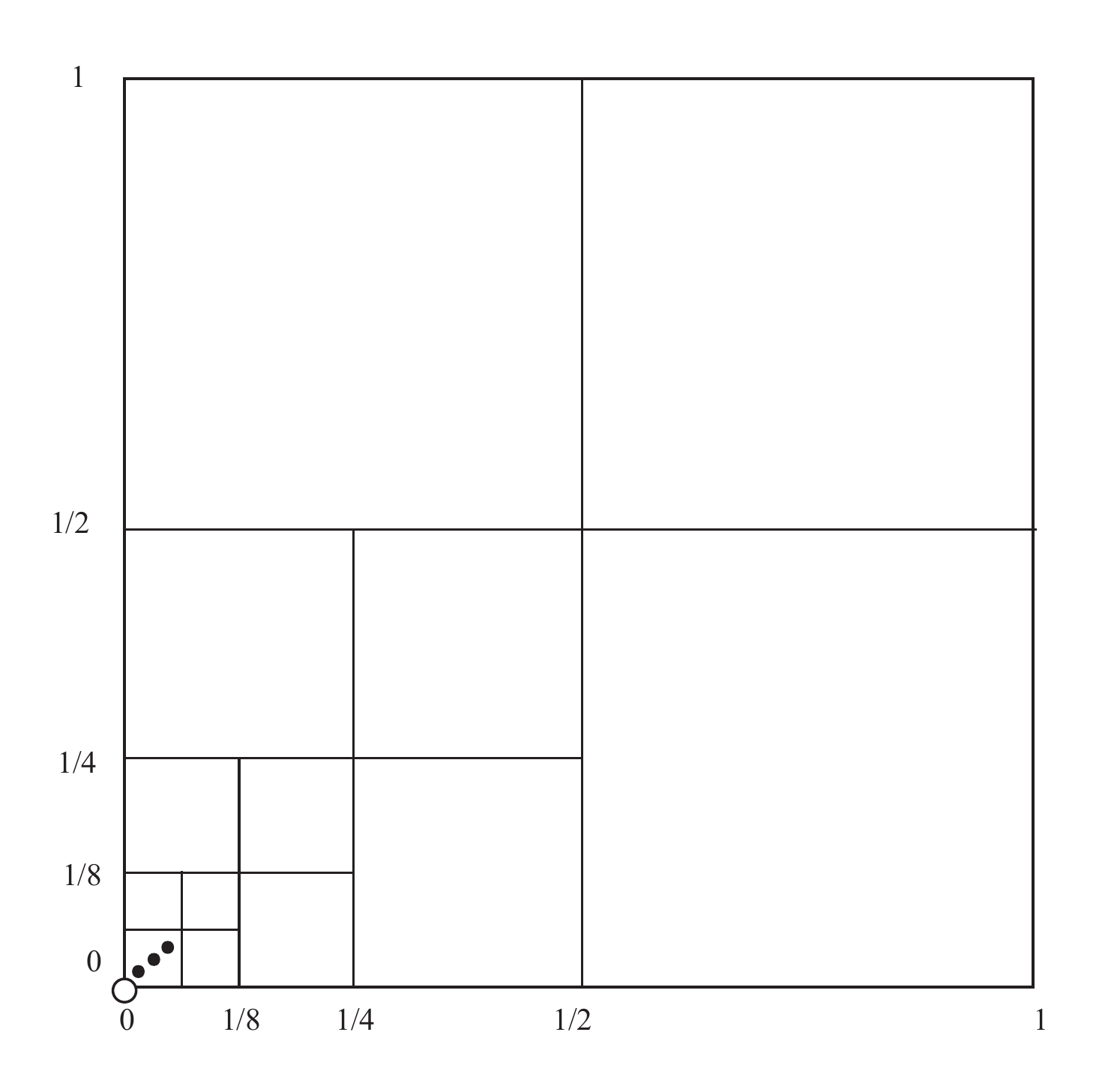}
\caption{A CW decomposition of a unit cube minus one corner.}
\label{fig:F3.8}
\end{figure}
\end{example}

Following the example, if two 1-cells that are in $Y$ meet at a 0-cell that is not in $Y$, and these two 1-cells are not in the same vertical line, we may apply an analogous cell decomposition as for $[0,1]\times[0,1] \setminus \{(0,0)\}$. 
Suppose these two 1-cells are on the same vertical line, we return to the exceptional case as discussed before in which we employ semicircles to further divide up the basic region. 

There is one more situation that is actually `troublesome'. Previously we've seen in Proposition~\ref{prop:F3.2}, if $Y$ is in the form of (\ref{eqn:F3.1}), it is impossible to have a basic region whose vertical sides consist of more than one 1-cells such that not all 1-cells are simultaneously inside $Y$ (see Figure~\ref{fig:F3.11}). However, if $Y$ is a union of at least two sets in the form of (\ref{eqn:F3.1}), this circumstance might not longer be true. 
For example, let $Y = {x > 0} \ \cup \  {y > 0}$, and consider the cell decomposition of $[-1, 1] \times [-1,1] \ \cap \ Y$. (That is to say, $R = 1$ in the proof of Proposition~\ref{prop:F3.3}.) Then the basic region on the right has its left side consisting of two 1-cells one of which is in $Y$ while the other is not. Our previous technique of overlaying cell decompositions with respect to the 1-cells on the boundary of a basic region fails here. 
\begin{figure}[ht]
\includegraphics[width=7cm]{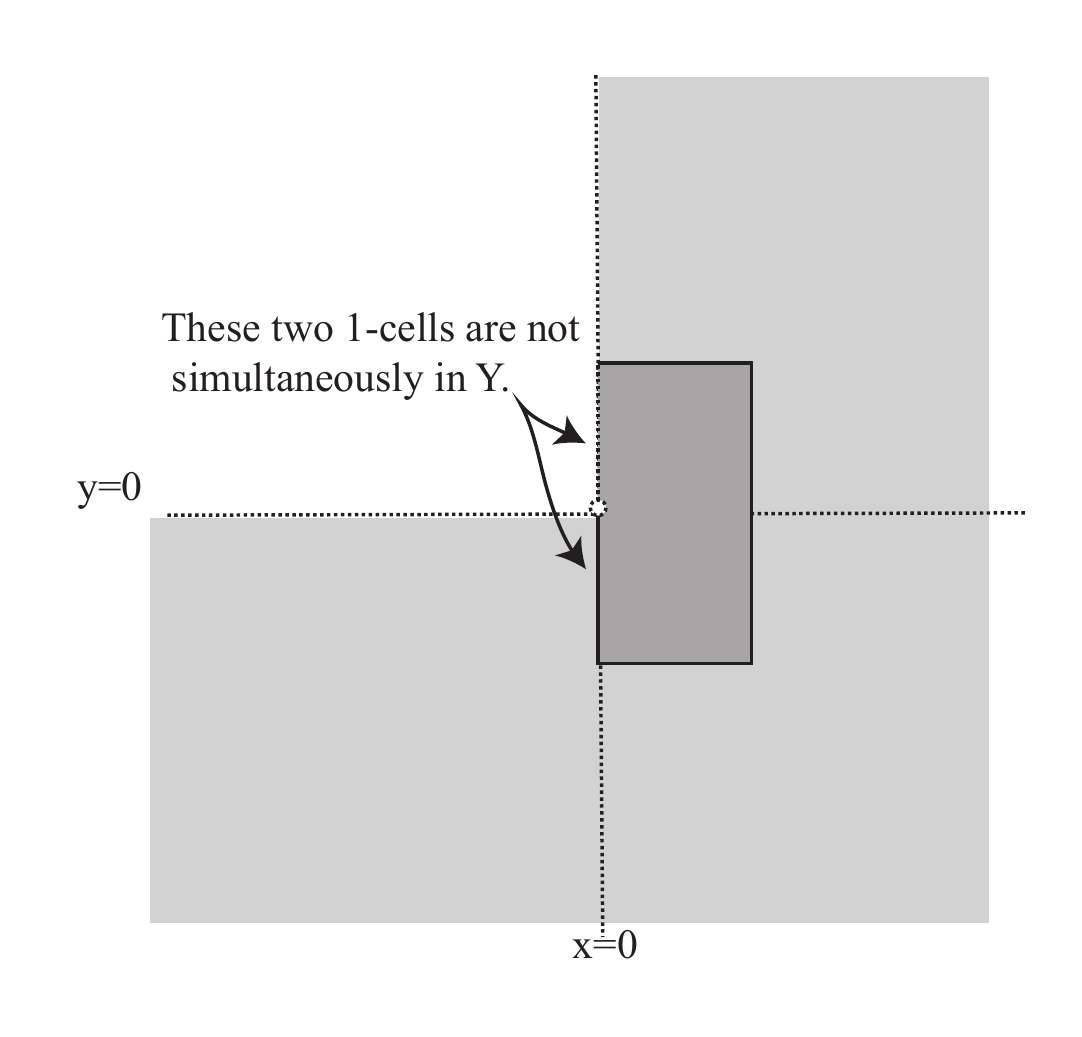}
\caption{If $Y = {x > 0} \cup {y > 0}$, there is a basic region with its left side consisting of two 1-cells one of which is in $Y$ while the other is not.}
\label{fig:F3.11}
\end{figure}
In this situation, let us insert an additional horizontal 1-cell at the 0-cell that connects two 1-cells one of which is in $Y$ while the other is not. Suppose this open line interval is completely contained in the basic region (that is, intersecting the top and bottom graphs at most at a point on the other side), then we can divide our original basic region into two basic regions, eliminating this `troublesome' case. 

However, it is very likely that such a horizontal line interval has a nontrivial intersection with the top or bottom graph of the region before even reaching the other side (see Figure~\ref{fig:F3.12}). Based on our construction, the top and bottom graphs belong to one of the following four types: 
\begin{enumerate}
\item{strictly increasing, convex upward}; 
\item{strictly decreasing, convex upward}; 
\item{convex downward}; 
\item{linear}.
\end{enumerate}
This is because, previously only inflection points and local minimum points were considered for the 0-cells, and it was sufficient. But here let us also include the local maximum points as 0-cells in the case of convex downward (Dr. Hardt communicated this to me). This will greatly simplify the argument, since the top and bottom graphs now belong to one of the following three types:
\begin{enumerate}
\item{strictly increasing};
\item{strictly decreasing};
\item{constant}.
\end{enumerate}

If the horizontal line interval intersects the top or bottom graph at a point which is not on the other side, either the top graph is strictly decreasing, or the bottom graph is strictly increasing. Furthermore, since the two graphs don't intersect except possibly at the other side, the horizontal 1-cell at the connecting 0-cell intersects only one of the top and bottoms graphs at exactly one point. 
\begin{figure}[ht]
\includegraphics[width=5cm]{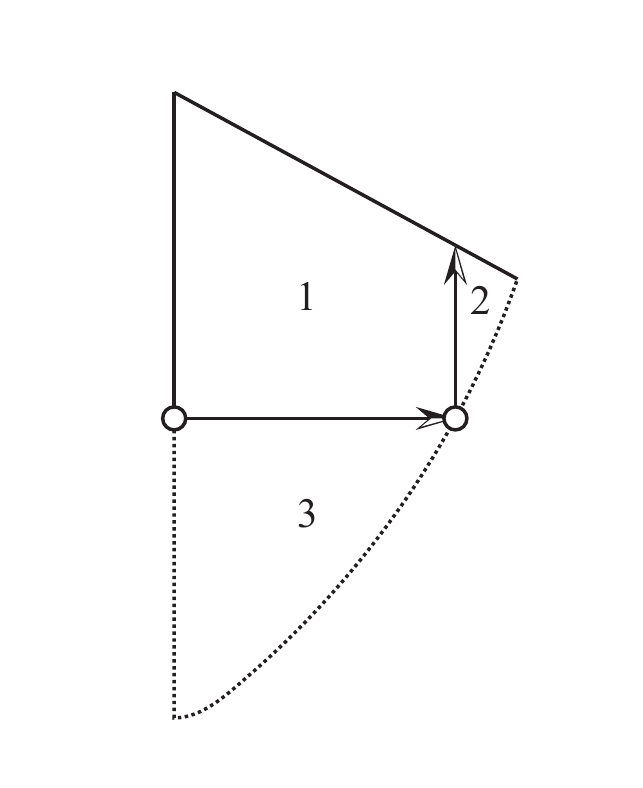}
\caption{Divide a basic region into three basic regions using a horizontal 1-cell at the troublesome 0-cell and a vertical 1-cell at the intersection point.}
\label{fig:F3.12}
\end{figure}
Therefore, we may introduce another vertical line interval at the intersection point, so that the original basic region can be divided into three basic regions, each of which is no longer `troublesome' (see Figure~\ref{fig:F3.12}).
If more than one such connecting 0-cells were present, we may perform the above procedure consecutively for each one of them. 

As a conclusion, there exists a cell decomposition $\mathcal{A}$ of $Y$ such that $(Y, \mathcal{A})$ is a CW complex and satisfies the finiteness property. Indeed, the technique of dividing in the `troublesome' case makes sure that each cell is still semi-algebraic and the local finiteness property also holds. 
\end{proof}

\begin{remark}
In Proposition~\ref{prop:F3.2}, we've shown how to use semicircles to decompose a basic region if it has a side with a removed 0-cell between two 1-cells. In fact, we may also divide the basic region using the technique above. 
\end{remark}

\section{The intersection of an algebraic set and an open semi-algebraic set}

Suppose $f_1, \ldots, f_k$, $g_1, \ldots, g_m$ are nonzero real-valued polynomials in two variables, and assume they are irreducible and distinct. Consider the following semi-algebraic set $Z$: 
\begin{equation}\label{eqn:F4.1}
Z = \{f_1(x, y)=0, \ldots, f_k(x, y)=0\} \cap \{g_1(x, y)>0, \ldots, g_m(x, y)>0\}. 
\end{equation}
If $k \geq 2$, the first $f_i$'s determine at most finitely many points, thus $Z$ is a finite set. Now let's assume that $k=1$. Define $B$ as follows:
\begin{equation}\label{eqn:F4.2}
B = \{(x, y) \in \mathbb{R}^2: f_1(x, y) = 0, g_j(x, y)=0 \text{ for some } j \in \{1, \ldots, m\}\}.
\end{equation}
Since $B$ is at most a finite set, we can include each point in $B$ as a 0-cell to the cell decomposition of $\{f_1(x, y)=0\}$ as guaranteed by  Corollary~\ref{cor:F2.2}. It follows that $V(f_1) \setminus B$ is a union of 0- and 1-cells. However, this is not yet a cell decomposition for some 1-cells' endpoints might be in $B$. To fix this problem, for each of these 1-cells which contain at least one endpoint in the set $B$, we replace it with infinitely many 1-cells. The idea can be best illustrated by looking at the following two examples. 
\begin{example}
The open unit interval $(0,1)$ may be decomposed as below: 
\begin{eqnarray}\label{eqn:F4.3}
\text{ 0-cells}: & \frac{1}{2^n}, \ \ 1-\frac{1}{2^n} \\ 
\text{ 1-cells}: & (\frac{1}{2^{n+1}}, \frac{1}{2^{n}}), \ \ (1-\frac{1}{2^n}, 1-\frac{1}{2^{n+1}}), \text{ where } n\geq 1.  \nonumber
 \end{eqnarray}
 And the half-closed interval $(0, 1]$ may be decomposed as below:
  \begin{eqnarray}\label{eqn:F4.4}
\text{ 0-cells}: & \frac{1}{2^n}, \ \ 1 \\ 
\text{ 1-cells}: & (\frac{1}{2^{n+1}}, \frac{1}{2^{n}}), \ \ (\frac{1}{2}, 1), \text{ where } n\geq 1.  \nonumber
 \end{eqnarray}
 \end{example}
 So if an endpoint of a 1-cell is removed, we consider a sequence of 0-cells converging to the endpoint; and the parts between consecutive 0-cells determine the infinitely many 1-cells. 
 As a result, $V(f_1) \setminus B$ becomes a cell complex. Furthermore, it is locally finite, thus is a CW complex. 
 
 For the semi-algebraic set $Z$, it consists of the cells in $V(f_1) \setminus B$ that are also in $\{g_1(x, y) > 0$, $\ldots$, $g_m(x, y) > 0\}$. We note that if a 1-cell in $V(f_1) \setminus B$ is in $Z$, then its two endpoints must be in $Z$ too, otherwise they are in $B$, which is a contradiction. Thus $Z$ is a subcomplex of $V(f_1)\setminus B$. The finiteness property for $Z$ is automatic. Let's summarize our result in the following proposition. 
 
 \begin{proposition}\label{prop:F4.1}
 Let $Z$ be a semi-algebraic set as defined in (\ref{eqn:F4.1}), there exists a cell decomposition $\mathcal{A}$ for $Z$ such that $(Z, \mathcal{A})$ satisfies the finiteness property which is also a CW complex. 
 \end{proposition}
 
 Let's take a finite union of these sets and see what happens. 
 \begin{proposition}\label{prop:F4.2}
 Let $Z$ be a finite union of semi-algebraic sets as defined in (\ref{eqn:F4.1}), then there exists a cell decomposition $\mathcal{A}$ for $Z$ such that $(Z, \mathcal{A})$ satisfies the finiteness property which is also a CW complex. 
 \end{proposition}
 
 \begin{proof}
 Suppose $Z$ is defined as follows:
 \begin{eqnarray*}
 Z & = & \{f_1=0, \ldots, f_k=0, g_1>0, \ldots, g_m>0\} \cup \{f'_1=0, \ldots, f'_{k'}=0, g'_1>0, \ldots, g'_{m'}>0\} \\
   & & \cup \ldots \cup \{f''_1=0, \ldots, f''_{k''}=0, g''_1>0, \ldots, g''_{m''}>0\},
 \end{eqnarray*}
 where the union is finite. We know that when $k \geq 2$, there are at most finitely many points. Therefore the above expression can be reduced to the following:
 \begin{eqnarray}\label{eqn:F4.5}
 Z & = & \{f_1=0, g_1>0, \ldots, g_m>0\} \cup \{f'_1=0, g'_1>0, \ldots, g'_{m'}>0\} \\
   & &\cup \ldots \cup \{f''_1=0, g''_1>0, \ldots, g''_{m''}>0\} \cup \{\text{ finitely many points }\}. \nonumber
 \end{eqnarray}
 
Consider the finite union $X$ of the affine algebraic sets $V(f_1)$, $V(f'_1)$, \ldots, $V(f''_1)$, together with the finitely many points in (\ref{eqn:F4.5}), then  $X$ has a CW decomposition. Next consider the union $\tilde{B}$ of sets in the form of (\ref{eqn:F4.2}):  
 \begin{eqnarray*}
\tilde{B} &=& \{f_1= 0, g_j =0, 1 \leq j \leq m\}  \cup \{f'_1= 0, g'_{j'}=0, 1 \leq j' \leq m' \} \\
& &\cup \ldots \cup \{f''_1 = 0, g''_{j''} =0, 1 \leq j'' \leq m''\}.
\end{eqnarray*} 
It follows that $\tilde{B}$ is a finite set. Include the points in $\tilde{B}$ as 0-cells to the cell decomposition of $X$. Call it $\mathcal{A}_0$. Then remove the 0-cells from $X$ that are in $\tilde{B}$, and are not in the union $Z$. Fix these 1-cells whose endpoints are removed as before. Thus $X \setminus B \cap Z^c $ gets a new CW decomposition. Call it $\mathcal{A}$. 
Pick these 0- and 1-cells in $\mathcal{A}$ that are contained $Z$, it turns out that $Z$ is a subcomplex of $\mathcal{A}$. Indeed, 
it suffices to show that if a 1-cell is contained in $Z$, its endpoints are contained in $Z$ too. Let's go back to the first cell decomposition $\mathcal{A}_0$. 
Given a 1-cell $e \in \mathcal{A}_0$, without loss of generality, we may assume that $e$ comes from $\{f_1=0\}$. If its endpoints are not in $\tilde{B}$, we pass it directly to $\mathcal{A}$. It follows that if $e$ is in $Z$, then its two endpoints are also in $Z$. 
Now let us suppose that at least one endpoint of $e$ is in $\tilde{B}$, say $p$. 
On the one hand, if $p$ is not in $Z$, then $p$ is removed from $\mathcal{A}$, and $e$ is replaced by infinitely many 1-cells, each of which returns to the previous case. On the other hand, if $p$ is in $Z$, then we keep $p$ in $\mathcal{A}$ and the end of $e$ connecting to $p$ remains intact. Therefore if $e$ is in $Z$, then $p$ is automatically in $Z$. Repeating the same argument for 1-cells coming from $\{f'_1=0\}$, \ldots, $\{f''_1=0\}$, yields the desired conclusion that $Z$ is a subcomplex. The finiteness property for $Z$ is easy to check. 
\end{proof}


\section{general case}
In general, an arbitrary semi-algebraic set $X$ in the plane can be described as: 
\begin{equation*}
X = \bigcup_{i=1}^{I} \bigcap_{j=1}^{J} \{(x, y) \in \mathbb{R}^2: f_{i, j}(x,y) = 0, g_{i,j}(x,y ) > 0\}, 
\end{equation*}
where $f_{i,j}$, $g_{i,j}$ are nonzero real-valued polynomials in two variables. We see that $X$ is a finite union of sets in the form obtained by taking the intersection of an algebraic set (i.e. $\{f_1(x, y)= 0, \ldots, f_k(x,y)=0\}$) with an open semi-algebraic set (i.e. $\{g_1(x, y) > 0, \ldots, g_m(x, y)>0\}$). What's more, we may assume that the $f_i$, $g_j$ are irreducible because of the following observations: 
\begin{eqnarray*}
\{f(x,y) \cdot \tilde{f}(x,y) = 0\} &= &\{f(x,y) = 0\} \bigcup \{\tilde{f}(x,y) = 0\}; \\
\{g(x,y) \cdot \tilde{g}(x,y) > 0\} &= &\big(\{g(x,y)>0\}\cap \{\tilde{g}(x,y)>0\} \big) \bigcup \\
						& & \big( \{-g(x,y) > 0 \} \cap \{-\tilde{g}(x,y) >0\} \big). 
\end{eqnarray*}

From previous results, we've known how to construct a CW decomposition with the finiteness property for each of the following three types of semi-algebraic sets: 
\begin{eqnarray*}
(I) & \{f_1(x, y) = 0, \ldots, f_k(x,y)=0\}, \text{ (Theorem~\ref{thm:F2.2}) }\\
(II) & \{g_1(x, y) > 0, \ldots, g_m(x, y)>0\}, \text{ (Theorem~\ref{thm:F3.4}) } \\
(III) & \{f_1(x, y) = 0, \ldots, f_k(x,y) = 0\} \cap \{g_1(x, y) > 0, \ldots, g_m(x, y)>0\}, \text{ (Proposition~\ref{prop:F4.1})}. 
\end{eqnarray*}

Now we are ready to take their finite unions. First, let's begin with finite unions of the same type. We've discussed them already in previous sections. Namely, Corollary~\ref{cor:F2.3} for a finite union of type (I); Theorem~\ref{thm:F3.4} for a finite union of type (II); and Proposition~\ref{prop:F4.2} for a finite union of type (III). 

Next, let's take a finite union of exactly two different types: (I) + (II), (I) + (III), and (II) + (III). 

\begin{lemma}[I + II]\label{lem:F5.1}
Suppose $W$ is a finite union of sets in the form of (I) and (II), then $W$ has a cell decomposition $\mathcal{A}$ that satisfies the finiteness property. However, $(W, \mathcal{A})$ is not necessarily a CW complex. 
\end{lemma}

\begin{proof}
From hypothesis, $W$ is in the following form: 
\begin{equation}\label{eqn:F5.1}
W = \{f_1=0\} \cup \ldots \cup \{f_k=0\} \cup \{g_1> 0, \ldots, g_m>0\} \cup \ldots \cup \{g'_1>0, \ldots, g'_{m'}>0\}.
\end{equation}
Let $X$ be defined as below:
\begin{eqnarray}\label{eqn:F5.2}
 X &=& \{f_1=0\} \cup \ldots \cup \{f_k=0\} \cup \{g_1= 0\} \cup \ldots \cup \{g_m=0\} \\ 
 & & \cup \ldots \cup \{g'_1=0\} \cup \ldots \cup \{g'_{m'}=0\}. \nonumber
 \end{eqnarray}
As usual, $X$ has a CW decomposition; moreover, we can divide $W$ up into basic regions. The cell decomposition is almost the same as in Propositions~\ref{prop:F3.2} and \ref{prop:F3.3}, with only one exception.
First, assume a basic (open) region (without boundary) is contained in $W$. 
In Proposition~\ref{prop:F3.3}, we see that if a side does not lie in the semi-algebraic set $Y$, then its two endpoints also do not lie in $Y$, which is essential for the local finiteness property on the boundary of a basic region. However, such a nice observation fails for $W$ here, in particular at corner points or at isolated removed points on vertical sides. A simple counterexample is $W = \{x^2+y^2=0\} \cup \{y > 0\}$ (see Figure~\ref{fig:F14}). 
\begin{figure}[ht]
\includegraphics[width=8cm]{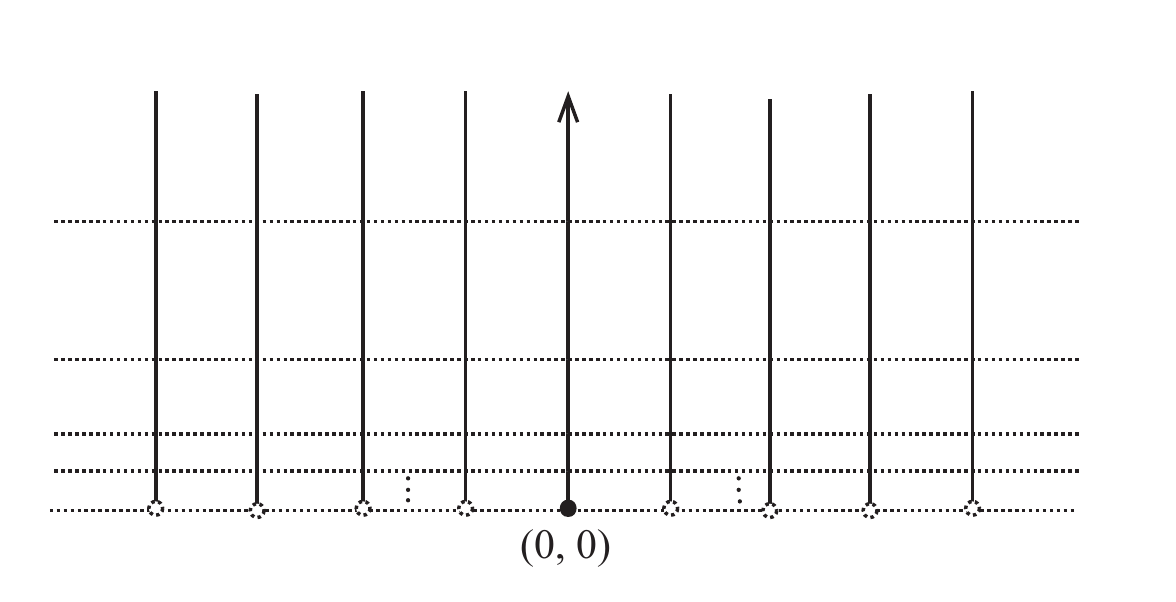}
\caption{For $W = \{x^2+y^2=0\} \cup \{y > 0\}$, the cell decomposition is not locally finite at the origin.}
\label{fig:F14}
\end{figure}
Thus it is possible to have 1-cell on the boundary of a basic region which is not in $W$ but either of whose endpoints is in $W$. If such a situation happens, for example, at a corner or at a removed isolated point on a vertical side, we must include this endpoint as a 0-cell, causing the local finiteness property to fail definitely. 

Second, assume a basic region is not contained in $W$. Then let's look at the cells on its boundary. A 1-cell is in $W$ if and only if it is contained in $V(f_1) \cup \ldots \cup V(f_k)$, which is closed. Therefore the two endpoints of the 1-cell are both contained in $W$, yielding a cell decomposition. 

Therefore, a cell decomposition exists for $W$ for which the local finiteness property might fail. However, the finiteness property can be checked to still hold. 
\end{proof}

\begin{lemma}[I + III]\label{lem:F5.2}
Suppose $W$ is a finite union of sets in the form of (I) and (III), then $W$ has a cell decomposition $\mathcal{A}$ that satisfies the finiteness property. Moreover, $(W, \mathcal{A})$ is a CW complex. 
\end{lemma}

\begin{proof}
The proof is analogous to that of Proposition~\ref{prop:F4.2}. More precisely,  
$W$ takes the following form by hypothesis:
\begin{eqnarray*}
W &=& \{f_1=0\} \cup \ldots \cup \{f_k=0\} \cup \{\tilde{f}_1 = 0, g_1 > 0, \ldots, g_m>0\} \\
& & \cup \ldots \cup \{\tilde{f}'_1=0, g'_1>0, \ldots, g'_{m'}>0\} \cup \{\text{ finitely many points }\}.
\end{eqnarray*}
Let $X$ be defined as below:
\begin{equation*}
X = \{f_1=0\} \cup \ldots \cup \{f_k=0\} \cup \{\tilde{f}_1 = 0\} \cup \ldots \cup \{\tilde{f}'_1=0\} \cup \{\text{ finitely many points }\},
\end{equation*}
which has a CW decomposition. Furthermore, add the following set $\tilde{B}$ as 0-cells to the cell decomposition. 
\begin{equation*}
\tilde{B} = \{\tilde{f}_1= 0, g_j =0, 1 \leq j \leq m\}  \cup \ldots \cup \{\tilde{f}'_1 = 0, g'_{j'} =0, 1 \leq j' \leq m'\}.
\end{equation*}  
Then we remove the 0-cells that are in $\tilde{B}$ and are not in $W$. Fixing these 1-cells whose endpoints are removed as in (\ref{eqn:F4.3}) and (\ref{eqn:F4.4}), and selecting those carried by $W$ yields a CW complex for $W$, which is a subcomplex of $X \setminus \tilde{B} \cap W^c$. It remains to check that if a 1-cell $e$ is contained in $W$, then its endpoints are also contained in $W$. In our construction, $e$ is carried entirely by one of the following affine algebraic sets:
\begin{equation*}
\{f_1=0\}, \ldots, \{f_k=0\}, \{\tilde{f}_1=0\}, \ldots, \{\tilde{f}'_1=0\}. 
\end{equation*} 

There are two cases. {\it Case~1:} $e$ is carried by $\{f_1=0\} \cup \ldots \cup \{f_k=0\}$, then its endpoints are automatic in $W$ by closedness. {\it Case~2:} $e$ is carried by $\{\tilde{f}_1=0\} \cup \ldots \cup \{\tilde{f}'_1=0\}$.  The endpoints are also contained in $W$, and the proof is analogous to that in Proposition~\ref{prop:F4.2}.
\end{proof}

\begin{lemma}[II + III]\label{lem:F5.3}
Suppose $W$ is a finite union of sets in the form of (II) and (III), then $W$ has a cell decomposition $\mathcal{A}$ that satisfies the finiteness property. However, $(W, \mathcal{A})$ might not be a CW complex. 
\end{lemma}

\begin{proof}
The hypothesis says that $W$ can be written in the following form:
\begin{eqnarray*}
W &=& \{g_1 >0, \ldots, g_m>0\} \cup \ldots \cup \{g'_1>0, \ldots, g'_{m'}>0\} \cup \\
& & \{f_1=0, \tilde{g}_1>0, \ldots, \tilde{g}_{\tilde{m}}>0\} \cup \ldots \cup \{f'_1=0, \tilde{g}'_1>0, \ldots, \tilde{g}'_{\tilde{m}'}>0\} \\
& & \cup \ \{\text{ finitely many points }\}.
\end{eqnarray*}
By Lemma~\ref{lem:F5.1}, there exists a cell decomposition with the finiteness property for the following set $W'$: 
\begin{eqnarray*}
W' &=& \{g_1 >0, \ldots, g_m>0\} \cup \ldots \cup \{g'_1>0, \ldots, g'_{m'}>0\} \\
& & \cup \ \{f_1=0\} \cup \ldots \cup \{f'_1=0\} \cup \{\text{ finitely many points }\},
\end{eqnarray*}
where every point $(x_0, y_0)$ in the plane corresponds to a real algebraic set such as $\{(x - x_0)^2+(y -y_0)^2 = 0 \}$. 
Furthermore, add the following set $\tilde{B}$ as 0-cells to the above cell decomposition (before dividing $W'$ into basic regions):
\begin{equation*}
\tilde{B} =  \{f_1= 0, \tilde{g}_j =0, 1 \leq j \leq \tilde{m}\}  \cup \ldots \cup \{f'_1 = 0, \tilde{g}'_{j'} =0, 1 \leq j' \leq \tilde{m}'\}.
\end{equation*}  
Then we need to remove the 0- and 1-cells that are not $W$, in particular these lying in the union $\{f_1=0\} \cup \ldots \cup \{f'_1=0\}$. Based on our construction, each of these cells is on the boundary of some basic region. There are two different cases. 

{\it Case 1: the basic region belongs to $W$.} Without loss of generality, we may assume that every side has no isolated 0-cells, otherwise dividing the region further by introducing horizontal and vertical line intervals as shown in Theorem~\ref{thm:F3.4}. It follows that every 0-cell is at the corner and every side is made up of only 1-cell. Thus removing a 0- or 1-cell won't affect the cell decomposition too much, except for some minor adjusts. In fact, we've seen all possible boundary conditions already. 

{\it Case 2: the basic region does not belong to $W$.} Removing a 1-cell won't affect anything. However, removing a 0-cell might cause a problem. Suppose $e$ is 1-cell adjacent to this 0-cell, and $e$ is in $W$. If $e$ is on the boundary of a basic region contained in $W$, we return to case 1. Otherwise, we need to fix this 1-cell (in particular, the half with the 0-cell as an endpoint) by replacing it with infinitely many smaller 1-cells and 0-cells according to (\ref{eqn:F4.4}). 

As a result, there exists a cell decomposition $\mathcal{A}$ for $W$. It is easy to check that $(W, \mathcal{A})$ satisfies the finiteness property. However, this cell decomposition is not necessarily locally finiteness for the same reason as shown in Lemma~\ref{lem:F5.1}. 
\end{proof}


Finally, we are ready to take a finite union of all three different types: (I) + (II) + (III). 
\begin{lemma}[I + II + III]\label{lem:F5.4}
Suppose $W$ is a finite union of sets in the form of (I), (II) and (III), then $W$ has a cell decomposition $\mathcal{A}$ that satisfies the finiteness property. However, $(W, \mathcal{A})$ might not be a CW complex. 
\end{lemma}

\begin{proof}
We may write $W$ as follows: 
\begin{eqnarray*}
W &=& \{f_1=0\} \cup \ldots \cup \{f_k=0\} \cup \ \{g_1>0, \ldots, g_m>0\} \cup \ldots \cup \{g'_1>0, \ldots, g'_{m'}>0\} \\
& & \cup \ \{\tilde{f}_1=0, \tilde{g}_1>0, \ldots, \tilde{g}_{\tilde{m}}>0\} \cup \ldots \cup \{\tilde{f}'_1=0, \tilde{g}'_1>0, \ldots, \tilde{g}'_{\tilde{m}'}>0\}.
\end{eqnarray*}
Then the rest of the proof is similar to that of Lemma~\ref{lem:F5.3}. 
\end{proof}

\begin{theorem}\label{thm:F5.5}
Given any semi-algebraic set in the plane, it has a cell decomposition with the finiteness property.
\end{theorem}

\begin{proof}
Combine Lemmas~\ref{lem:F5.1}, \ref{lem:F5.2}, \ref{lem:F5.3}, and \ref{lem:F5.4}. 
\end{proof}

\section{Conclusion}
In this paper, we find a semi-algebraic stratification $\mathcal{A}$, in particular a cell decomposition, for any arbitrary semi-algebraic set $X$ in the plane. Moreover, $\mathcal{A}$ satisfies an analytic condition concerning geodesics. More precisely, suppose $A$, $B$ are two arbitrary points in $X$, and $\gamma$ is a piecewise $C^2$ curve from $A$ to $B$ lying entirely in $X$ such that its length is the shortest among all possible such curves. Then the intersection of $\gamma$ with every cell in $\mathcal{A}$ is either empty or consists of finitely many components, each of which is either a singleton or a geodesic line segment. 

Furthermore, when $X$ is in one of the following cases, $(X, \mathcal{A})$ turns out to be a CW complex, because the cell decomposition is locally finite. 
\begin{enumerate}
\item{$X$ is a finite union of sets in the form of $\{f_1=0, \ldots, f_k=0\}$;}
\item{$X$ is a finite union of sets in the form of $\{g_1>0, \ldots, g_m>0\}$;}
\item{$X$ is a finite union of sets in the form of $\{\tilde{f}_1=0, \ldots, \tilde{f}_k=0, \tilde{g}_1>0, \ldots, \tilde{g}_m>0\}$;}
\item{$X$ is a finite union of sets in the form of $\{f_1=0, \ldots, f_k=0\}$ and \\
$\{\tilde{f}_1=0, \ldots, \tilde{f}_k=0, \tilde{g}_1>0, \ldots, \tilde{g}_m>0\}$.}
\end{enumerate}

The future questions may concern higher dimensional semi-algebraic sets, or semi-analytic sets, or sub-analytic sets, or triangulations, or even more complicated analytical conditions such as Lipschitz conditions (that is to say, whether each 1-cell is the graph of a Lipschitz function). 
\bibliographystyle{amsplain}

\end{document}